\magnification 1200
\input amstex.tex
\input amsppt.sty
\voffset-4truemm
\hoffset4truemm

\author Nazarov and Tarasov \endauthor
\title Representations of Yangians \endtitle


\expandafter\ifx\csname maxim.def\endcsname\relax \else\endinput\fi
\expandafter\edef\csname maxim.def\endcsname{%
 \catcode`\noexpand\@=\the\catcode`\@\space}
\catcode`\@=11

\mathsurround 1.6pt

\def\hcor#1{\advance\hoffset by #1}
\def\vcor#1{\advance\voffset by #1}
\let\bls\baselineskip  \let\ignore\ignorespaces
\def\vsk#1>{\vskip#1\bls} \let\adv\advance 
\def\vv#1>{\vadjust{\vsk#1>}\ignore} \def\vvv#1>{\vadjust{\vskip#1}\ignore}
\def\vvn#1>{\vadjust{\nobreak\vsk#1>\nobreak}\ignore}
\def\vvvn#1>{\vadjust{\nobreak\vskip#1\nobreak}\ignore}
\def\emph#1{{\it #1\/}}

\let\vp\vphantom  
 \let\nt\noindent \let\cl\centerline
\def\nn#1>{\noalign{\vskip #1pt}} \def\NN#1>{\openup#1pt}
 
\let\Sum\sum \def\sum{\Sum\limits} 
\let\Prod\prod \def\prod{\Prod\limits} \let\Int\int \def\int{\Int\limits}

\let\=\m@th \def\&{.\kern.1em} \def\>{\!\;} \def\:{\!\!\;}

\ifx\plainfootnote\undefined \let\plainfootnote\footnote \fi
\expandafter\ifx\csname amsppt.sty\endcsname\relax
 
\else \fi

\newbox\s@ctb@x
\def\s@ct#1 #2\par{\removelastskip\vsk>
 \vtop{\bf\setbox\s@ctb@x\hbox{#1} \parindent\wd\s@ctb@x
 \ifdim\parindent>0pt\adv\parindent.5em\fi\item{#1}#2\strut}%
 \nointerlineskip\nobreak\vtop{\strut}\nobreak\vsk-.4>\nobreak}

\newbox\t@stb@x
\def\gadv{\global\advance} \def\gad#1{\gadv#1 1} 
\def\l@b@l#1#2{\def\n@@{\csname #2no\endcsname}%
 \if *#1\gad\n@@ \expandafter\xdef\csname @#1@#2@\endcsname{\the\Sno.\the\n@@}%
 \else\expandafter\ifx\csname @#1@#2@\endcsname\relax\gad\n@@
 \expandafter\xdef\csname @#1@#2@\endcsname{\the\Sno.\the\n@@}\fi\fi}
\def\l@bel#1#2{\l@b@l{#1}{#2}\?#1@#2?}
\def\?#1?{\csname @#1@\endcsname}
\def\[#1]{\def\n@xt@{\ifx\t@st *\def\n@xt####1{{\setbox\t@stb@x\hbox{\?#1@F?}%
 \ifnum\wd\t@stb@x=0 {\bf???}\else\?#1@F?\fi}}\else
 \def\n@xt{{\setbox\t@stb@x\hbox{\?#1@L?}\ifnum\wd\t@stb@x=0 {\bf???}\else
 \?#1@L?\fi}}\fi\n@xt}\futurelet\t@st\n@xt@}
\def\(#1){{\rm\setbox\t@stb@x\hbox{\?#1@F?}\ifnum\wd\t@stb@x=0 ({\bf???})\else
 (\?#1@F?)\fi}}
\def\dff{\expandafter\d@f} \def\d@f{\expandafter\def}
\def\edff{\expandafter\ed@f} \def\ed@f{\expandafter\edef}

\newcount\Sno \newcount\Lno \newcount\Fno
\def\Section#1{\gad\Sno\Fno=0\Lno=0\s@ct{\the\Sno.} {#1}\par} \let\Sect\Section
\def\section#1{\gad\Sno\Fno=0\Lno=0\s@ct{} {#1}\par} \let\sect\section
\def\l@F#1{\l@bel{#1}F} \def\<#1>{\l@b@l{#1}F} \def\l@L#1{\l@bel{#1}L}
\def\Tag#1{\tag\l@F{#1}} \def\Tagg#1{\tag"\llap{\rm(\l@F{#1})}"}
\def\Th#1{Theorem \l@L{#1}} \def\Lm#1{Lemma \l@L{#1}}
\def\Prop#1{Proposition \l@L{#1}}
\def\Cr#1{Corollary \l@L{#1}} \def\Cj#1{Conjecture \l@L{#1}}
 
\def\Proof#1.{\demo{\it Proof #1}}

\def\Par{\par\medskip} \def\setparindent{\edef\Parindent{\the\parindent}}
\def\Appendix{\Sno=64\let\p@r@\z@ 
\def\Section##1{\gad\Sno\Fno=0\Lno=0 \s@ct{} \hskip\p@r@ Appendix \char\the\Sno
 \if *##1\relax\else {.\enspace##1}\fi\par} \let\Sect\Section
\def\section##1{\gad\Sno\Fno=0\Lno=0 \s@ct{} \hskip\p@r@ Appendix%
 \if *##1\relax\else {.\enspace##1}\fi\par} \let\sect\section
\def\l@b@l##1##2{\def\n@@{\csname ##2no\endcsname}%
 \if *##1\gad\n@@
\expandafter\xdef\csname @##1@##2@\endcsname{\char\the\Sno.\the\n@@}%
\else\expandafter\ifx\csname @##1@##2@\endcsname\relax\gad\n@@
 \expandafter\xdef\csname @##1@##2@\endcsname{\char\the\Sno.\the\n@@}\fi\fi}}

\let\logo@\relax
\let\m@k@h@@d\makeheadline \let\m@k@f@@t\makefootline
\def\makeheadline{\ifnum\pageno=1\headline={\hfil}\fi\m@k@h@@d}
\def\makefootline{\ifnum\pageno=1\footline={\hfil}\fi\m@k@f@@t}


\def\E(#1){\mathop{\hbox{\rm End}\,}(#1)} 
\def\id{\hbox{\rm id}}  \def\for{\hbox{for \,}}


\let\al\alpha

 \let\eps\varepsilon \let\epsilon\eps
\let\ka\kappa
\let\la\lambda \let\La\Lambda
\let\si\sigma  
 \let\phi\varphi
\let\om\omega \let\Om\Omega

\def\T{\Bbb T}
\def\Z{\Bbb Z}


\def\deri{\text{\rm d}}
\def\imin{{\text{\bf\~\i}}}
\def\jmin{{\text{\bf\~\j}}}

\def\al{\alpha}
\def\AN{{\operatorname{A}(\glN)}}
\def\ANg{{\operatorname{A}_g(\glN)}}

\def\be{\beta}
\def\bi{{\bold i}}
\def\bj{{\bold j}}
\def\bk{{\bold k}}
\def\bl{{\bold l}}
\def\bm{{\bold m}}
\def\bn{{\bold n}}
\def\bx{{\boxed{\phantom{\square}}\kern-.4pt}}

\def\CC{{\Bbb C}}

\def\ddd{\leaders\hbox to.6em{\hss.\hss}\hfil}
\def\de{\delta}
\def\De{\Delta}
\def\deg{{\operatorname{deg}\ts}}

\def\End{{\operatorname{End}}}
\def\EndCN{{\End(\CC^N)}}
\def\enddemos{{$\quad\square$\enddemo}}

\def\id{{\operatorname{id}}}
\def\io{\iota}

\def\ga{\gamma}
\def\gz{\vcenter{
 \hbox to 150pt{$\la_{N1}\ \ \ \la_{N2}\ddd\la_{NN}$}
 \hbox to 145pt{$\kern10pt\la_{N-1,1}\ \la_{N-1,2}\ddd\la_{N-1,N-1}$}
 \hbox to 94pt{$\kern22pt\ddots\hfil\ndots$}
 \hbox to 84pt{$\kern37pt\la_{21}\hfil\la_{22}$}
 \hbox{$\kern50pt\la_{11}$}}}
 \def\ndots{\mathinner{\mkern1mu\raise1pt\hbox{.}\mkern2mu\raise4pt
 \hbox{.}\mkern2mu\raise7pt\vbox{\kern7pt\hbox{.}}\mkern1mu}}
\def\gl{{\frak{gl}}}
\def\glM{{\frak{gl}_M}}
\def\glMN{{\frak{gl}_{M+N}}}
\def\glN{{\frak{gl}_N}}

\def\ka{\kappa}

\def\la{\lambda}
\def\La{\Lambda}
\def\lc{{,\ts\ldots,\ts}}
\def\lm{{\la,\mu}}

\def\mi{{\raise.5pt\hbox{-}}}
\def\mv{{\kern55pt}}
\def\mw{{\kern84pt}}

\def\of{{\om_f}}
\def\ot{\otimes}
\def\om{\omega}
\def\Om{\Omega}

\def\ph{\varphi}
\def\ps{\psi}

\def\S{{\Cal S}}
\def\si{\sigma}

\def\Slm{{\S_\lm}}
\def\slN{{\frak{sl}_N}}
\def\Sn{{\S_\la}}

\def\T{{\Cal T}}
\def\Tlm{{\T_\lm}}
\def\tm{{\ts-}}
\def\Tn{{\T_\la}}
\def\tp{{\ts+}}
\def\ts{\thinspace}
\def\tz{{\ts0}}

\def\UM{{\operatorname{U}(\glM)}}
\def\UMN{{\operatorname{U}(\glMN)}}
\def\UN{{\operatorname{U}(\glN)}}
\def\UNh{{\operatorname{U}_q(\widehat{\gl}_N)}}
\def\UNx{{\operatorname{U}(\glN[u])}}
\def\Up{\Upsilon}

\def\X{{\Cal X}}
\def\Xn{{\operatorname{X}^{(n)}}}

\def\YM{{\operatorname{Y}(\glM)}}
\def\YMN{{\operatorname{Y}(\glMN)}}
\def\YN{{\operatorname{Y}(\glN)}}
\def\YsN{{\operatorname{Y}(\slN)}}

\def\Z{{\Cal Z}}
\def\ZM{{\operatorname{Z}(\glM)}}
\def\ZN{{\operatorname{Z}(\glN)}}
\def\ZZ{{\Bbb Z}}

\document

\font\bigf=cmr10 scaled 1200
\cl{\bigf Representations of Yangians with Gelfand-Zetlin bases}
\bigskip
\bigskip
\cl{By {\it Maxim Nazarov\,} at York \,and
{\it Vitaly Tarasov\,} at St.\,Petersburg}
\bigskip
\bigskip
\centerline{\hbox to 4cm{\hrulefill}}
\bigskip
\bigskip
\nt
We study a certain family of finite-dimensional modules
over the Yangian $\YN$. The Yangian $\YN$ is a canonical deformation of the
universal enveloping algebra $\UNx$
in the class of Hopf algebras [D1].
The classification of the irreducible finite-dimensional $\YN$-modules
has been obtained in
[D2] by generalizing the results of the second author [\ts T1,T2\ts]. These
modules
are parametrized by all pairs formed by a sequence of monic polynomials
$P_1(u)\lc P_{N-1}(u)\in\CC[u]$ and a  power series
$f(u)\in1+u^{-1}\ts\CC[[u^{-1}]]$.   

The algebra $\YN$ comes equipped with a distinguished maximal commutative
subalgebra $\AN$. This subalgebra is generated by the centres of all
algebras in the chain
$$
\operatorname{Y}(\frak{gl}_1)
\subset
\operatorname{Y}(\frak{gl}_2)
\subset\ldots\subset
\YN.
$$
Continuing [C], we study
finite-dimensional $\YN$-modules with a semisimple action of the
subalgebra $\AN$. We will call these modules tame.
We provide a characterization
of the irreducible tame modules
in terms of the polynomials $P_1(u)\lc P_{N-1}(u)$;
see~Theorem 4.1.

It turns out that the spectrum of the action of the subalgebra $\AN$ in
every
irreducible tame module is simple. The eigenbases of $\AN$ in the latter
modules
are called Gelfand-Zetlin bases. In Section 3 we provide explicit formulas
for the action of the generators of the algebra $\YN$ introduced in [D2]
on the vectors of these bases. Moreover, each of these bases contains
the vector $\xi_0$
called singular.
For every vector $\xi$ of the Gelfand-Zetlin basis we point out an element
$b\in\YN$ 
such that $b\cdot\xi_0=\xi$.

The algebra $\YN$ admits a homomorphism onto the universal enveloping
alge\-bra
$\UN$. Denote this homomorphism by $\pi_N$.
The image of the centre of $\YN$ with respect to $\pi_N$
coincides with the centre $\ZN$ of the algebra $\UN$. So any  
finite-dimensional $\YN$-module obtained via the
homomorphism $\pi_N$ is tame. The results of Section 3 can be regarded
as a generalization of the classical formulas of [GZ].
For more details on the connection between the formulas of [GZ]
and the Yangian $\YN$
see our publication [NT]. It has inspired the publication~[M].

For each $M=0,1,2,\ldots$ there is a homomorphism of
$\YN$~to the commutant in $\UMN$ of the subalgebra $\UM$. The image of this
homomorphism along with the centre $\ZM$ generates the commutant [O].
For any dominant integral weights $\la$ and $\mu$ of the Lie algebras
$\glMN$ and $\glM$ denote by $V_\lm$ 
the subspace in the irreducible $\glMN$-module of the highest weight
$\la$ formed by all singular vectors with respect to $\glM$ of the weight
$\mu$.
The latter homomorphism makes~$V_\lm$ into an irreducible tame $\YN$-module.

Furthermore, for any $h\in\CC$ there is an automorphism of the algebra $\YN$
preserving the subalgebra $\AN$.
Applying this automorphism to an irreducible
finite-dimensional $\YN$-module $V$ amounts to changing $u$ by $u+h$
in the respective polynomials $P_1(u)\lc P_{N-1}(u)$ and the series $f(u)$.
Denote the resulting $\YN$-module by $V(h)$. We prove that
up to the choice of the series $f(u)$
every irreducible tame $\YN$-module splits into a tensor product 
of modules of the form $V_\lm(h)$; see Theorem 4.1. 

Consider the $\YN$-module $V(h)$ 
where the module $V$ is  
obtained from the the fundamental $\gl_N$-module $\CC^N$ 
via the homomorphism $\pi_N$. 
The above stated splitting property of any irreducible tame $\YN$-module
was established
in [C], Theorem~3.9 under two extra conditions.
First, the spectrum of the action of the algebra $\AN$ in this module was
assumed to be simple. Second, this module was assumed to be a subquotient
of the module $V(h)\ot\ldots\ot V(h^\prime)$ where
the number of tensor factors is less than $N$. It was also conjectured in
[C] that both these conditions could be removed. So we confirm this
conjecture in the present article.

The Yangian $\YN$ can be viewed [D3] as a degeneration
of the quantum univ\-er\-sal
enveloping algebra $\UNh$. All the results presented in this article have
their quantum counterparts. We will give them in a forthcoming
publication.
Note that
Theorem 4.11 from [GRV] is the
counterpart of our formulas for the action of $\YN$ in the module
$V(h)\ot\ldots\ot V(h^\prime)$ where
$V=\CC^N$ as above while the parameters $h\lc h^\prime\in\CC$
are in general position;
see Theorem 2.9 and Corollary~3.9 here.

It is our pleasure to thank 
Ivan Cherednik and Grigori Olshanski for numerous creative conversations. 
The first author was supported by the University of Wales Research
Fellowship.
The second author should like to thank the Mathematics De\-partment of the
University of Wales at Swansea for hospitality.

\bigskip\bigskip
\Section{Preliminaries on Yangians}
\nt
In this section we gather several known facts about the {\it Yangian}
of the Lie algebra $\glN$. This is a complex associative
unital algebra $\YN$ with the
countable
set of generators $T^{(s)}_{ij}$ where $s=1,2,\ts\ldots$ and
$i,j=1,\ts\dots\ts,N$. The de\-fin\-ing relations in the algebra $\YN$ are
$$
[\ts T_{ij}^{(r)},T_{kl}^{(s+1)}\ts]-
[\ts T_{ij}^{(r+1)},T_{kl}^{(s)}\ts]=
T_{kj}^{(r)}\ts T_{il}^{(s)}-T_{kj}^{(s)}\ts T_{il}^{(r)};\qquad 
r,s=0,1,2,\ldots\ts
\Tag{1.1}
$$
where $T_{ij}^{(0)}=\de_{ij}\cdot1$.
We will also
use the following matrix form of these relations.
\Par
Let $E_{ij}\in\EndCN$ be the standard matrix units.
Introduce two formal variables $u,v$ and consider the
{\it Yang $R$-matrix}
$$
R(u,v)=(u-v)\cdot\id+\sum_{i,j=1}^N\ts E_{ij}\ot E_{ji}
\in\EndCN^{\ot\ts2}\ts[u,v].  
$$
Introduce the formal power series in $u^{-1}$
$$
T_{ij}(u)=
T_{ij}^{(0)}+T_{ij}^{(1)}\ts u^{-1}+T_{ij}^{(2)}\ts u^{-2}+\ldots
$$
and combine all these series into the single element
$$
T(u)=\sum^N_{i,j=1}\ts E_{ij}\ot T_{ij}(u)
\in\EndCN\ot\YN\ts[[u^{-1}]].
$$

For any associative unital algebra $\operatorname{X}$ denote 
by $\io_s$ its embedding into a finite tensor product
$\operatorname{X}^{\ot\ts n}$
as the $s$-th tensor factor:
$$
\io_s(x)=1^{\ot\ts (s-1)}\ot x\ot1^{\ot\ts(n-s)},\quad
x\in\operatorname{X};\qquad s=1,\dots,n.
$$
Introduce also the formal power series with the coefficients 
in $\EndCN^{\ot2}\ot\YN$
$$
T_1(u)=\io_1\ot\id\bigl(T(u)\bigr),
\quad
T_2(v)=\io_2\ot\id\bigl(T(v)\bigr),
\quad
R_{12}(u,v)=R(u,v)\ot1.
$$
Then the relations \(1.1) can be rewritten as
$$
R_{12}(u,v)\ts T_1(u)\ts T_2(v)=T_2(v)\ts T_1(u)\ts R_{12}(u,v)
\Tag{1.2}
$$ 

The relations \(1.2) imply that
for any formal series $f(u)\in 1+u^{-1}\ts\CC[[u^{-1}]]$ 
the assignement of the generating series
$T_{ij}(u)\mapsto f(u)\cdot T_{ij}(u)$
determines an automorphism of the algebra $\YN$. We will denote by $\of$ this
automorphism.
The element $T(u)$ of the algebra $\EndCN\ot\YN\ts[[u^{-1}]]$
is invertible; denote
$$
T(u)^{-1}=\widetilde T(u)=\sum^N_{i,j=1}\ts E_{ij}\ot\widetilde{T}_{ij}(u).
\Tag{1.3}
$$
Then the relations \(1.2) along with the equality
$
R(u,v)\ts R(-u,-v)=1-(u-v)^2
$
imply that the assignment $T_{ij}(u)\mapsto\widetilde{T}_{ij}(-u)$ 
determines an automorphism of the algebra $\YN$.
We will denote  by $\si_N$ this automorphism; it is clearly involutive. 
We will also make use of the automorphism $\tau_N$ determined by the
assignement
$$
T_{ij}(u)\mapsto\widetilde{T}_{N-j+1,N-i+1}(u).
$$

\proclaim{Lemma 1.1}
Suppose that $i\neq l$ and $k\neq j$. Then the coefficients of the series
$T_{ij}(u)$ commute with those of $\widetilde T_{kl}(u)$.
\endproclaim

\demo{Proof}
Multiplying \(1.2) on the left and on the right by $T_2(v)^{-1}$
we obtain the equality
$$
T_1(u)\ts R_{12}(u,v)\ts T_2(v)^{-1}=T_2(v)^{-1}\ts R_{12}(u,v)\ts T_1(u).
$$ 
Since $i\neq l$ and $k\neq j$ by multiplying the latter equality by
$E_{ii}\ot E_{kk}\ot1$ on the left
and by $E_{jj}\ot E_{ll}\ot1$ on the right we obtain that
$$
E_{ij}\ot E_{kl}\ot
\bigl[\ts T_{ij}(u)\ts,\widetilde T_{kl}(-v)\ts\bigr]=0
\quad\square
$$
\enddemo

\nt
Due to (1.1)
every sequence of pairwise distinct indices $\bk=(k_1,\dots,k_M)$ where
$1\leqslant k_i\leqslant M+N$,
determines an embedding of algebras
$$
\ph^{\vp1}_\bk:\YM\to\YMN:\ts T_{ij}(u)\mapsto T_{k_i k_j}(u).
\Tag{1.4}
$$
We will also make use of the embedding of the same algebras
$
\ps^{\vp1}_\bk=\si_{M+N}\ts\ph^{\vp1}_\bk\ts\si_M.
$
The embedding $\ph^{\vp1}_\bk$ with $\bk=(1,\dots,M)$ will be called
{\it standard}. We will identify the algebra $\YM$ and its image in $\YMN$
with respect to the standard embedding.

Now let $\bk=(k_1,\dots,k_M)$ and $\bl=(k_1,\dots,k_N)$
be two sequences such that
$$\{\ts k_1,\dots,k_M\ts\}
\sqcup
\{\ts l_1,\dots,l_N\ts\}=
\{1,\dots,M+N\ts\}.
\Tag{1.5}
$$

\proclaim{Corollary 1.2}
The images of the embeddings $\ph^{\vp1}_\bk$ and
$\ps^{\vp1}_\bl$ in $\YMN$ commute.
\endproclaim

\demo{Proof}
The image of the mapping $\ps^{\vp1}_\bl$ coincides
with that of $\si_{M+N}\ts\ph^{\vp1}_\bl$ and we have
$
\si_{M+N}\ts\ph^{\vp1}_\bl\bigl(T_{mn}(u)\bigr)=
\widetilde T_{l_m l_n}(u).
$
Due to Lemma 1.1 the coefficients of the series
$T_{k_i k_j}(u)$ commute with those of $\widetilde T_{l_m l_n}(u)$
for all possible $i,j,m,n$
\enddemos

\nt
Consider the elements $E_{ij}$ as generators of the Lie algebra $\glN$.
The algebra $\YN$ contains the universal enveloping elgebra $\UN$
as a subalgebra: due to \(1.1) the assignment $E_{ji}\mapsto T_{ij}^{(1)}$
defines the embedding. Moreover, there is a homomorphism
$$
\pi_N:\YN\to\UN:\ts T_{ij}(u)\mapsto \de_{ij}+E_{ji}\ts u^{-1}.
\Tag{1.505}
$$
Note that this homomorphism is by definition identical on the subalgebra
$\UN$. By \(1.3) the automorphism $\si_N$ of the algebra $\YN$ is also
identical on $\UN$. 

Let $\ZM$ denote the center of the universal enveloping algebra $\UM$. 
By the definition \(1.4) the image
$\ph^{\vp1}_\bk\bigl(\UM\bigr)$ is contained in $\UMN$. For the proof of
the following proposition see [O].

\proclaim{Proposition 1.3}
The centralizer of the subalgebra $\ph^{\vp1}_\bk\bigl(\UM\bigr)$
in $\UMN$ coincides with
$$\ph^{\vp1}_\bk\bigl(\ZM\bigr)\cdot
\pi_{M+N}\ts\ps^{\vp1}_\bl\bigl(\YN\bigr).$$
\endproclaim

\nt
We will now describe an alternative set of generators for the algebra $\YN$.
These generators were introduced in [D2], and we will call them
{Drinfeld generators}. Let 
$\bi=(i_1,\dots,i_m)$ and $\bj=(j_1,\dots,j_m)$ be two
sequences of indices
such that 
$$
1\leqslant i_1<\ldots<i_m\leqslant N
\quad\text{and}\quad
1\leqslant j_1<\ldots<j_m\leqslant N.
\Tag{1.55}
$$
Consider the sum over all permutations $g$ of $1,2,\dots,m$
$$
Q^{\vp1}_{\bi\bj}(u)=
\sum_g\ts
T_{i_1 j_{g(1)}}(u)\ts T_{i_2 j_{g(2)}}(u-1)\ \ldots\ T_{i_m j_{g(m)}}(u-m+1)
\cdot\operatorname{sign}g\ts;
$$
here  the series
in $(u-1)^{-1},\ts\dots\ts,(u-m+1)^{-1}$ should be
re-expanded in $u^{-1}$.
Denote
$$
A_m(u)=Q^{\vp1}_{\bi\bi}(u)\quad\for\bi=(1,\ldots,m).
\Tag{1.505}
$$
We will set $A_0(u)=1$.
The series $A_N(u)$ is called the
{\it quantum determinant} for the algebra $\YN$. 
The following proposition is well known; a detailed proof can be found in
[MNO], Section 2.

\proclaim{Proposition 1.4}
The coefficients at $u^{-1},\ts u^{-2},\ts\dots$ of the series $A_N(u)$
are free generators for the centre of the algebra $\YN$.
\endproclaim 

\nt
Consider the ascending chain of algebras
$$
\operatorname{Y}(\frak{gl}_1)
\subset
\operatorname{Y}(\frak{gl}_2)
\subset\ldots\subset
\YN
\Tag{1.6}
$$
determined by the standard embeddings.  
Denote by $\AN$ the subalgebra of $\YN$ generated by the centres of all the
algebras in \(1.6); this subalgebra is commutative. By Proposition 1.4 the
coefficients of the series $A_1(u),\dots,A_N(u)$ generate the subalgebra
$\AN$.
The subalgebra $\AN$ is maximal commutative in $\YN$.
This fact is contained in [C], Theorem 2.2. However, we will not use
this fact in the present article.
Now for every $m=1\lc N-1$ denote
$$
B_m(u)=Q^{\vp1}_{\bi\bj}(u),\quad
C_m(u)=Q^{\vp1}_{\bj\bi}(u),\quad
D_m(u)=Q^{\vp1}_{\bj\bj}(u)
\Tag{1.605}
$$
where $\bi=(1,\ldots,m)$ and $\bj=(1,\ldots,m-1,m+1)$.
The coefficients of the series 
$$
A_m(u),\quad B_m(u),\quad C_m(u);\qquad m=1\lc N-1
$$
along with those of $A_N(u)$ generate the algebra $\YN$; see
[D2], Example.
We will call
the coefficients of all these series the {\it Drinfeld generators} for $\YN$.

\demo{Remark}
Consider the fixed point subalgebra in $\YN$ with respect to all
automorphisms
of the form $\of$. This algebra is called the {\it Yangian} for the
simple Lie
algebra $\slN$ and denoted by $\YsN$.  
The coefficients of all the series   
$$
A_{m-1}(u-1)\ts A_m^{-1}(u-1)\ts A_m^{-1}(u)\ts A_{m+1}(u),
\quad
A_m^{-1}(u)\ts B_m(u),
\quad
C_m(u)\ts A_m^{-1}(u)
$$
with $m=1\lc N-1$
are generators of the algebra $\YsN$; see again [D2], Example.
The algebra $\YN$ is isomorphic
to the tenzor product of its centre and the algebra $\YsN$; see
[MNO], Section 2 for
the proof of this claim.    
\enddemo

\nt
Now let $\bi$ and $\bj$ be any two sequences of indices satisfying 
the condition \(1.55). Introduce the increasing sequences 
{\bf\u\i}$=(i_{m+1}\lc i_N)$ and
{\bf\u\j}$=(j_{m+1}\lc j_N)$
such that
$$
\{\ts i_1\lc i_m\ts\}
\sqcup
\{\ts i_{m+1}\lc i_N\ts\}=
\{\ts j_1\lc j_m\ts\}
\sqcup
\{\ts j_{m+1}\lc j_N\ts\}=
\{1\lc N\ts\}.
$$
Let $\epsilon$ be the sign of the permutation
$(i_1\lc i_n)\mapsto (j_1\lc j_n)$.
Introduce also the sequences
$$
\align
&\imin=(N-i_N+1\lc N-i_{m+1}+1),
\\
&\jmin=(N-j_N+1\lc N-j_{m+1}+1).
\endalign
$$
\proclaim{Lemma 1.5}
In the algebra $\YN[[u^{-1}]]$ we have the equalities
$$
\align
\si_N\bigl(Q_{\bi\bj}(u)\bigr)
&=\epsilon\cdot
Q_{\text{{\bf\u \j}\ts{\bf\u \i}\ts}}(-\ts1-u)\ts/\ts A_N(m-1-u),
\Tag{1.6333}
\\
\tau_N\bigl(Q_{\bi\bj}(u)\bigr)
&=\epsilon\cdot
Q_{\imin\ts\jmin\ts}(u-m)\ts/\ts A_N(u).
\Tag{1.6666}
\endalign
$$
\endproclaim

\demo{Proof}
Let $H\in\EndCN^{\ot m}$ denote the antisymmetrization map
normalized so that $H^2=m!\ts H\ts$. Then introduce the map 
$$
H_m=\io_1\ldots\io_m(H)\in\EndCN^{\ot m}.
$$
Introduce also the
formal power series with the coefficients in $\EndCN^{\ot N}\ot\YN$
$$
T_1(u)=\io_1\ot\id\bigl(T(u)\bigr)
\ ,\ \ldots\ ,\ 
T_N(v)=\io_N\ot\id\bigl(T(u)\bigr).
$$
Then we have
$$
\align
&T_1(u)\ts T_2(u-1)\ \ldots\ T_m(u-m+1)\cdot(H_m\ot1)=
\Tag{1.6999}
\\
&(H_m\ot1)
\cdot
T_m(u-m+1)\ \ldots\ T_2(u-1)\ts T_1(u)\ts;
\endalign
$$
see [MNO], Section 2 for the proof of this equality. In particular,
when $m=N$
this equality implies that
$$
T_1(u)\ts T_2(u-1)\ \ldots\ T_N(u-N+1)\cdot(H_N\ot1)=
H_N\ot A_N(u).
$$
Therefore we have
$$
\align
T_{m+1}(u-m)\ts T_{m+2}(u-m-1)
\ \ldots\ T_N(u-N+1)
\cdot(H_N\ot1)&=
\\
\widetilde T_m(u-m+1)
\ \ldots \ \widetilde T_2(u-1)\ts\widetilde T_1(u)
\cdot\bigl(H_N\ot A_N(u)\bigr)&\ts.
\Tag{1.6777}
\endalign
$$
By multiplying the latter equality on the left by 
$I_N\ot1$ where $I_N$ equals
$$
E_{i_mi_m}\ot\ldots\ot E_{i_1i_1}\ot
E_{j_{m+1}j_{m+1}}\ot\ldots\ot E_{j_Nj_N}
$$
we obtain that
$$
\epsilon\cdot Q_{\text{{\bf\u \j}\ts{\bf\u \i}\ts}}(u-m)=
\si_N\bigl(Q_{\bi\bj}(m-1-u)\bigr)\ts A_N(u). 
$$
Replacing here the variable $u$ by $m-1-u$ we obtain \(1.6333).
The equality \(1.6999) also implies
$$
\align
&
\widetilde T_m(u-m+1)\ \ldots\ \widetilde T_2(u-1)\ts \widetilde T_1(u)
\cdot(H_m\ot1)=
\\
&
(H_m\ot1)\cdot
\widetilde T_1(u)\ts\widetilde T_2(u-1)\ \ldots\ \widetilde T_m(u-m+1).
\endalign
$$
Therefore by multiplying the equality \(1.6777)  on the left by 
$J_N\ot1$ where $J_N$ equals
$$
\align
E_{N-j_1+1,N-j_1+1}\ot&\ldots\ot E_{N-j_m+1,N-j_m+1}\ot
\\
E_{N-i_N+1,N-i_N+1}\ot&\ldots\ot E_{N-i_{m+1}+1,N-i_{m+1}+1}
\endalign
$$
we obtain that
$$
\epsilon\cdot
Q_{\imin\ts\jmin\ts}(u-m)=
\tau_N\bigl(Q_{\bi\bj}(u)\bigr)\ts A_N(u)\quad\square 
$$
\enddemo  

\nt
By making use of Lemma 1.1 and Proposition 1.4 along with the equality
\(1.6333) we obtain
the following corollary to Lemma 1.5.

\proclaim{Corollary 1.6}
\!The coefficients of the series $Q_{\bi\bj}(u)$ commute with those of
$T_{i_kj_l}(u)\!$ for all $k,l=1\lc m$.
\endproclaim

\nt
Using the definitions \(1.505) and \(1.605)
we obtain one more corollary to Lemma 1.5.

\proclaim{Corollary 1.7}
In the algebra $\YN[[u^{-1}]]$ we have the equalities
$$
\align
\tau_N\bigl(A_m(u)\bigr)
&=\ts\phantom{-}
A_{N-m}(u-m)\ts/\ts A_N(u);
\qquad
m=1\lc N;
\Tag{1.655555}
\\
\qquad
\tau_N\bigl(B_m(u)\bigr)
&=-\ts
B_{N-m}(u-m)\ts/\ts A_N(u);
\qquad
m=1\lc N-1;
\\
\tau_N\bigl(C_m(u)\bigr)
&=-\ts
C_{N-m}(u-m)\ts/\ts A_N(u);
\qquad
m=1\lc N-1;
\Tag{1.666666}
\\
\tau_N\bigl(D_m(u)\bigr)
&=\ts\phantom{-}
D_{N-m}(u-m)\ts/\ts A_N(u);
\qquad
m=1\lc N.
\endalign
$$
\endproclaim

\demo{Remark}
The automorphism $\tau_N$ of the algebra $\YN$ is not involutive. Lemma 1.5
shows that
$$
\tau_N^{\ts2}\bigl(T_{ij}(u)\bigr)=T_{ij}(u-N)\cdot A_N(u)/A_N(u-1).
$$
For more comments on the automorphism $\tau_N^{\ts2}$ of $\YN$
 see [MNO], Section 5. 
\enddemo

\nt
Now let an arbitrary $M=0,1,2,\dots$ be fixed. Denote
$$
\align
\text{{\bf\B\i}}&=(1\lc M,M+i_1\lc M+i_m),
\\
\text{{\bf\B\j}}&=(1\lc M,M+j_1\lc M+j_m).
\endalign
$$

\proclaim{Corollary 1.8}
In the algebra $\YMN[[u^{-1}]]$ we have the equality
$$
Q_{\text{{\bf\B\i}\ts{\bf\B\j}}}(u)=
\ps^{\vp1}_\bl\bigl(Q_{\bi\bj}(u)\bigr)\cdot A_M(u-m)
\quad\for\quad 
\bl=(M+1\lc M+N).
$$
\endproclaim

\demo{Proof}
Introduce the sequences
$$
\bm=(M+i_{m+1}\lc M+i_N)
\quad\text{and}\quad
\bn=(M+j_{m+1}\lc M+j_N)\ts.
$$
By making use of Lemma 1.5 we then obtain the equalities
$$
\gather
\ps^{\vp1}_\bl
\bigl(Q_{\bi\bj}(u)\bigr)=
\epsilon\cdot
\si_{M+N}\ts\ph^{\vp1}_\bl\ts\si_N
\bigl(Q_{\bi\bj}(u)\bigr)=
\\
\epsilon\cdot
\si_{M+N}\ts\ph^{\vp1}_\bl
\bigl(Q_{\text{{\bf\u \j}\ts{\bf\u \i}\ts}}(-\ts1-u)\ts/\ts A_N(m-1-u)\bigr)
\\
=\si_{M+N}\bigl(Q_{\bn\bm}(-1-u)\cdot Q_{\bl\bl}(m-1-u)^{-1}\bigr) 
\\
=\bigl(Q_{\text{{\bf\B\i}\ts{\bf\B\j}}\ts}(u)\ts/\ts A_{M+N}(u+N-m)\bigr)
\cdot
\bigl(A_M(u-m)\ts/\ts A_{M+N}(u+N-m)\bigr)^{-1}
\\
=Q_{\text{{\bf\B\i}\ts{\bf\B\j}}\ts}(u)
\cdot
A_M(u-m)^{-1}\qquad\square
\endgather
$$
\enddemo

\demo{Remark}
By making use of Corollary 1.8
we obtain for each $m=1\lc N-1$ the following four equalities in the algebra
$\YN[[u^{-1}]]$:
$$
\align
A_m(u)&=A_{m-1}(u-1)\cdot\ps_{m,m+1}\bigl(\ts T_{11}(u)\bigr),
\\
B_m(u)&=A_{m-1}(u-1)\cdot\ps_{m,m+1}\bigl(\ts T_{12}(u)\bigr),
\\
C_m(u)&=A_{m-1}(u-1)\cdot\ps_{m,m+1}\bigl(\ts T_{21}(u)\bigr),
\\
D_m(u)&=A_{m-1}(u-1)\cdot\ps_{m,m+1}\bigl(\ts T_{22}(u)\bigr).
\endalign
$$
These equalities can be regarded as an alternative description of
the Drinfeld generators for the algebra $\YN$.
\enddemo

\proclaim{Proposition 1.9}
We have the commutation relations in $\YN[[u^{-1},v^{-1}]]$ 
$$
\align
[A_m(u),B_n(v)]&=0\quad\text{if}\quad m\neq n,
\Tag{1.7}
\\
[C_m(u),B_n(v)]&=0\quad\text{if}\quad m\neq n,
\Tag{1.8}
\\
[B_m(u),B_n(v)]&=0\quad\text{if}\quad|m-n|\neq1,
\Tag{1.9}
\\
\noalign{\vskip3pt}
(u-v)\cdot[A_m(u),B_m(v)]&=B_m(u)\ts A_m(v)-B_m(v)\ts A_m(u),
\Tag{1.10}
\\
(u-v)\cdot[C_m(u),B_m(v)]&=D_m(u)\ts A_m(v)-D_m(v)\ts A_m(u).
\Tag{1.11}
\endalign
$$
\endproclaim

\demo{Proof}
Let $\bm$ and $\bn$ be any subsequences of $\bi$ and $\bj$
respectively. Suppose that $\bm$ and $\bn$ are of the same length.
Then by Corollary 1.6 the coefficients of the
series $Q_{\bm\bn}(u)$  commute with those of $Q_{\bi\bj}(u)$.
This proves \(1.7) to \(1.9). For the proof of the relations \(1.10) to
\(1.11)
see [NT], Section 1
\enddemos

\nt
The proof of the next proposition is also contained in [NT], Section 1.

\proclaim{Proposition 1.10}
The following relation holds in the algebra $\YN[[u^{-1}]]$:
$$
C_m(u)\ts B_m(u-1)=D_m(u)\ts A_m(u-1)-A_{m+1}(u)\ts A_{m-1}(u-1).
\Tag{1.12}
$$
\endproclaim

\nt
There is a natural Hopf algebra structure on $\YN$ [D1]. The comultiplication
$\YN\to\YN^{\ot2}$ is defined by the assignement of the generating series
$$
T_{ij}(u)\mapsto\sum_{k=1}^N\ts
T_{ik}(u)\ot T_{kj}(u).
\Tag{1.13}
$$
Here and in what follows we take the tensor product over $\CC[[u^{-1}]]$
of the elements from the algebra $\YN[[u^{-1}]]$.

We will now consider the images of the Drinfeld generators for
the algebra $\YN$ with respect to the comultiplication
$$
\De^{(n)}:\YN\to\YN^{\ot n}.
\Tag{1.135}
$$
Let $\bi$ and $\bj$ be any two sequences of indices satisfying 
the condition \(1.55).

\proclaim{Proposition 1.11}
We have the equality
$$
\De^{(n)}\bigl(Q_{\bi\bj}(u)\bigr)=
\sum_{\bk^{(1)},\bk^{(2)}\lc\bk^{(n-1)}}\ts
Q_{\bi\bk^{(1)}}(u)\ot
Q_{\bk^{(1)}\bk^{(2)}}(u)\ot\ldots\ot
Q_{\bk^{(n-1)}\bj}(u).
$$
where 
$\bk^{(1)},\bk^{(2)}\lc\bk^{(n-1)}$ run through all the increasing sequences
of the inte\-gers $1\lc N$ of the length $m$.
\endproclaim

\demo{Proof}
We employ arguments similar to those used in the proof of Lemma 1.5.
Let $H\in\EndCN^{\ot m}$ be the antisymmetrization map inroduced therein.
Introduce the
formal power series with the coefficients in $\EndCN^{\ot m}\ot\YN$
$$
T_k(u)=\io_k\ot\id\bigl(T(u)\bigr);
\qquad
k=1\lc m.
$$
Then we have the relation
$$
\align
m!\cdot T_1(u)\ts T_2(u-1)\ts\ldots\ts T_m(u-m+1)\cdot(H\ot1)&=
\Tag{1.15}
\\
(H\ot1)\cdot T_1(u)\ts T_2(u-1)\ts\ldots\ts T_m(u-m+1)\cdot(H\ot1)&\ts.
\endalign
$$
Put 
$$
I=E_{i_1i_1}\ot\ldots\ot E_{i_mj_m},
\ \ 
J=E_{j_1j_1}\ot\ldots\ot E_{j_mj_m},
\ \ 
K=E_{i_1j_1}\ot\ldots\ot E_{i_mj_m}.
$$
Then by the definition \(1.505)
we have the equality 
$$
K\ot Q_{\bi\bj}(u)=
(I\ot1)\cdot
T_1(u)\ts T_2(u-1)\ts\ldots\ts T_m(u-m+1)\cdot
(H J\ot1).
\Tag{1.14}
$$

Introduce the formal power series with coefficients in
$\EndCN^{\ot m}\ot\YN^{\ot n}$
$$
T_{k[s]}(u)=\io_k\ot\io_s\bigl(T(u)\bigr);
\qquad
k=1\lc m;
\quad
s=1\lc n
$$
where $\io_k$ and $\io_s$ are respectively the
embeddings $\EndCN\to\EndCN^{\ot m}$ and
$\YN\to\YN^{\ot n}$. For any $k=1\lc m$ by the
definition \(1.13) we then have
$$
\id\ot\De^{(n)}\bigl(T_k(u)\bigr)=
T_{k[1]}(u)\ts T_{k[2]}(u)\ts\ldots\ts T_{k[n]}(u).
$$
Therefore from the equality \(1.14) we obtain that  
$$
\align
K\ot\De^{(n)}\bigl(Q_{\bi\bj}(u)\bigr)=
(I\ot1)\cdot
T_{1[1]}(u)\ts T_{1[2]}(u)\ts\ldots\ts T_{1[n]}(u)&\times 
\\
T_{2[1]}(u-1)\ts T_{2[2]}(u-1)\ts\ldots\ts T_{2[n]}(u-1)&\times\ldots\times 
\\
T_{m[1]}(u-m+1)\ts T_{m[2]}(u-m+1)\ts\ldots\ts T_{m[n]}(u-m+1)&\cdot
(H J\ot1)
\endalign
$$
where $I\ot1$ and $H J\ot1$ are now elements of the product
$\EndCN^{\ot m}\ot\YN^{\ot n}$.

Observe that the coefficients of the series $T_{k[s]}$ commute with those
of $T_{l[r]}$ if $k\neq l$ and $s\neq r$. Thus the right hand side of the
latter equality coincides with
$$
\align
(I\ot1)\cdot
T_{1[1]}(u)\ts T_{2[1]}(u-1)\ts\ldots\ts T_{m[1]}(u-m+1)&\times 
\\
T_{1[2]}(u)\ts T_{2[2]}(u-1)\ts\ldots\ts T_{m[2]}(u-m+1)&\times\ldots\times
\\
T_{1[n]}(u)\ts T_{2[n]}(u-1)\ts\ldots\ts T_{m[n]}(u-m+1)&\times
(H J\ot1).
\endalign
$$
By applying the relation \(1.15) repeatedly we now obtain that
$$
\align
K\ot\De^{(n)}\bigl(Q_{\bi\bj}(u)\bigr)=(m!)^{n-1}\ &\times 
\\
(I\ot1)\cdot
T_{1[1]}(u)\ts T_{2[1]}(u-1)\ts\ldots\ts T_{m[1]}(u-m+1)&\cdot(H\ot1)\times 
\\
T_{1[2]}(u)\ts T_{2[2]}(u-1)\ts\ldots\ts T_{m[2]}(u-m+1)&\cdot(H\ot1)
\times\ldots\times 
\\
T_{1[n]}(u)\ts T_{2[n]}(u-m+1)\ts\ldots\ts T_{m[n]}(u-m+1)&\cdot
(H J\ot1).
\endalign
$$
Proposition 1.11 follows from this equality and from the relation \(1.15)
\enddemos
 
\proclaim{Corollary 1.12}
We have the equality
$\De^{(n)}\bigl(A_N(u)\bigr)=\bigl(A_N(u)\bigr)^{\ot n}$.
\endproclaim

\nt
Note that $\De^{(n)}\bigl(A_m(u)\bigr)\neq\bigl(A_m(u)\bigr)^{\ot n}$
in general.
However, the next corollary will be sufficient for our purposes. 
Endow the algebra $\YN$ with the $\ZZ$-grading $\ts\deg$ determined by
$$
\deg T_{ij}^{(s)}=i-j;\qquad s=0,1,2,\dots\ts;
\Tag{1.999}
$$
see the relations \(1.1). We will extend this grading to the algebra
$\YN[[u^{-1}]]$ by assuming that $\deg u^{-1}=0$. The definitions \(1.505)
and \(1.605) then show that
$$
\deg A_m(u)=0,\quad \deg B_m(u)=-1,\quad \deg C_m(u)=1.
$$
The algebras $\YN^{\ot n}$ and $\YN^{\ot n}[[u^{-1}]]$ then
acquire grading by the group $\ZZ^n$. We will fix the lexicographical
ordering
on the group $\ZZ^n$.

\proclaim{Corollary 1.13}
For every $m=1\lc N$ we have the equality
$$
\De^{(n)}\bigl(A_m(u)\bigr)=\bigl(A_m(u)\bigr)^{\ot n}+
\ts\text{terms of the smaller degrees}.
$$
\endproclaim

\demo{Proof}
Let us apply Lemma 1.11 to $\bi=\bj=(1\lc m)$. Then 
$\deg Q_{\bi\bk^{(1)}}(u)<0$ unless $\bk^{(1)}=\bi$. Repeating this argument
we obtain the required statement
\enddemos

\nt
In Section 4 we will make use of the following observation.
For each permutation $g$ of $1,2\lc N$ denote by $\ANg$ the subalgebra in
$\YN$
generated by all the coefficients of the series $Q_{\bi\bi}(u)$ where
$$
\bi=\bigl(\ts1\lc N\ts\bigr)\setminus\ts g(m+1)\lc g(N)\ts\ts;
\qquad
m=1\lc N.
\Tag{1.999999999}
$$

\proclaim{Proposition 1.14}
Let $V$ be any finite-dimensional module of the algebra $\YN$.
The images in $\End(V)$ of all the subalgebras $\ANg$ are conjugated to each
other.
\endproclaim

\demo{Proof}
Denote by $\theta$ the action of the algebra $\YN$ in $V$.
Take the \text{embedding} of the algebra $\UN$ into $\YN$
defined by $E_{ji}\mapsto T_{ij}^{(1)}$. By the relations \(1.1) we have
$$
[\ts E_{ji}\ts,\ts T_{kl}^{(s)}\ts]
=
\de_{il}\cdot T_{kj}^{(s)}-\de_{kj}\cdot T_{il}^{(s)}\ts;
\qquad
s=1,2,\ldots\ts.
$$
Therefore the image of $\glN$ in $\End(V)$ contains an element $G$ such that
$$
\operatorname{exp}G
\cdot
\theta\bigl(\ts T_{ij}^{(s)}\ts\bigr)
=
\theta\bigl(\ts T_{g(i)\ts g(j)}^{\ts(s)}\ts\bigr)
\cdot  
\operatorname{exp}G\ts;
\qquad
s=1,2,\ldots\ts.\quad\ts\ts
$$
The latter equalities along with \(1.6999) imply that for each $m=1\lc N$
we have
$$
\operatorname{exp}G
\cdot
\theta\bigl(\ts A_m(u)\ts\bigr)
=
\theta\bigl(\ts Q_{\bi\bi}(u)\ts\bigr)
\cdot  
\operatorname{exp}G
$$
where the increasing sequence $\bi$ is defined in \(1.999999999)
\enddemos

\bigskip\bigskip
\Section{Modules with Gelfand-Zetlin bases}

\nt
In this section we will consider certain family of 
finite-dimensional modules over the algebra $\YN$.
Some of them will be irreducible and have a simple spectrum with respect to
the action of 
the maximal
commutative subalgebra $\AN$. The eigenbases of $\AN$ in the latter
modules will be called {Gelfand-Zetlin bases.}  We will determine the
eigenvalues of the coefficients of the series $A_1(u)\lc A_N(u)$
corresponding
to the vectors of these bases. The action of the remaining Drinfeld
generators
of the algebra $\YN$ on the vectors of these bases will be described in
Section 3.

Consider the matrix units $E_{ij}$ with $i,j=1\lc N$ as generators 
of the Lie algebra $\glN$. For every non-increasing sequence of integers
$\la=(\ts\la_1\lc\la_N\ts)$ denote by $V_\la$ the irreducible $\glN$-module
of the highest weight $\la\ts$. If $\xi\in V_\la$ is the highest
weight vector
then we have $E_{ii}\cdot\xi=\la_i\ts\xi$ and $E_{ij}\cdot\xi=0$ for $i<j$.

Denote by $\Tn$ the set of all arrays $\La$ with integral entries of the
form
$$
\gz
$$
where $\la_{Ni}=\la_i$ and $\la_{mi}\leqslant\la_i$ for all $i$ and $m$.
An array $\La\in\Tn$ is called a {\it Gelfand-Zetlin scheme} if
$
\la_{mi}\geqslant\la_{m-1,i}\geqslant\la_{m,i+1}
$
for all possible $m$ and $i$. Denote by $\Sn$ the subset in $\Tn$
consisting of the Gelfand-Zetlin schemes.

There is a canonical decomposition of the space $V_\la$ into the direct
sum of
one-dimensional subspaces associated with the acsending chain of subalgebras
$$
\operatorname{U}(\gl_1)
\subset
\operatorname{U}(\gl_2)
\subset\ldots\subset
\operatorname{U}(\glN)
\Tag{2.05}
$$
determined by the standard embeddings.
These subspaces are parametrized by the elements ot the set $\Sn$.
For each $m=1,2\lc N-1$
the subspace $V_\La\subset V_\la$ corresponding to the scheme $\La\in\Sn$
is contained in an irreducible
$\gl_m$-submodule of the highest weight $(\la_{m1},\la_{m2}\lc\la_{mm})$.
These conditions define the subspace $V_\La$ uniquely,
cf. [GZ]. The highest weight subspace in $V_\la$ corresponds to the scheme
$\La$ where $\la_{mi}=\la_i$ for every $m=1\lc N$.  

Consider the $\ZZ$-grading on the algebra $\UN$ such that the degree
of the element $E_{ij}$ equals $j-i$. It is the grading by the eigenvalues
of the adjoint action in $\UN$ of the element
$NE_{11}+(N-1)E_{22}+\ldots+E_{NN}$. Endow the space $V_\la$ with the
$\ZZ$-grading
by the eigenvalues of the action of the same element. Then 
the action of $\UN$ in the space $V_\la$ becomes a graded action.
Note that every subspace $V_\La$ in $V_\la$ is a weight subspace
and hence homogeneous; its degree equals
$$
\sum_{m=1}^N\ts\ts(N-m+1)
\left(\ts
\sum_{i=1}^m\ts\la_{mi}-\sum_{i=1}^{m-1}\ts\la_{m-1,i}
\ts\right)
=\sum_{m=1}^N\ts\sum_{i=1}^m\ts\la_{mi}\ts.
\Tag{2.07}
$$

Let us now make use of the homomorphism $\pi_N:\YN\to\UN$ defined by
\(1.505). 
Note that this homomorphism preserves the $\ZZ$-grading on the algebra $\YN$
introduced in the end of Section 1.
For each $m=1\lc N$ and any $\La\in\Tn$ introduce the rational function
$$
\al_{m\La}(u)=\prod_{i=1}^m\ts\frac
{\left(u+\la_{mi}-i+1\right)}{\left(u-i+1\right)}\ts;
\Tag{2.0733}
$$
note that
$\al_{m\La}(u)$ can be also regarded as a formal power series in $u^{-1}$.

\proclaim{Lemma 2.1}
The subspace $V_\La\subset V_\la$ is an eigenspace for the coefficients
of the
\nolinebreak
series $\pi_N\bigl(A_m(u)\bigr)$, the eigenvalues are the
respective coefficients of the series $\al_{m\La}(u)$.
\endproclaim

\demo{Proof}
By the definition \(1.505) we have
$$
\pi_N\bigl(A_m(u)\bigr)\in
\ \prod_{i=1}^m\ts\frac
{\left(u+E_{ii}-i+1\right)}{\left(u-i+1\right)}+
\operatorname{U}(\gl_m)[[u^{-1}]]\cdot\frak{n}_m
$$
where $\frak{n}_m$ is the subalgebra in $\gl_m$ spanned by the elements
$E_{ij}$ with $i<j$. 
Due to Proposition 1.4 the coefficients of the series
$\pi_N\bigl(A_m(u)\bigr)$
belong to the center of the algebra $\operatorname{U}(\gl_m)$ and act in any
irreducible $\gl_m$-submodule of $V_\la$ via scalars. 
Applying these coefficients to the highest weight
vector in an irreducible $\gl_m$-submodule of $V_\la$
we get the statement of Lemma 2.1 by the definition of the
subspace $V_\La$
\enddemos

\nt
Note that due to the conditions
$\la_{m1}\geqslant\la_{m2}\geqslant\ldots\geqslant\la_{mm}$ the  
scheme $\La\in\S_\la$ can be uniquely restored from the
collection of rational functions  
$\al_{1\La}(u)\lc\al_{N\La}(u)$.
Therefore the action in $V_\la$ of the commutative subalgebra
$\pi_N\bigl(\AN\bigr)$
in $\UN$ has a simple spectrum. 

\demo{Remark}
The subalgebra $\pi_N\bigl(\AN\bigr)$ in $\UN$
coincides with the subalgebra generated by the centres of all the 
algebras in \(2.05). The latter subalgebra in $\UN$ is maximal commutative
[C], Theorem 2.2.
\enddemo

\nt
Let $M$ run through the set $\{0,1,2,\dots\}$.
Fix the standard embedding of the algebra $\UM$ into $\UMN$.
For any pair of non-increasing sequences of integers
$$
\la=(\ts\la_1\lc\la_M,\la_{M+1}\lc\la_{M+N}\ts)
\quad\text{and}\quad
\mu=(\ts\mu_1\lc\mu_M\ts)
\Tag{2.0744}
$$
denote by $V_\lm$ the subspace in the $\glMN$-module $V_\la$ formed by all
the singular vectors with respect to $\glM$ of the weight $\mu$. 
Denote by $\Slm$ the subset in $\Sn$ formed by all the arrays 
$$
\La=(\ts\la_{li}\ts|\ts 1\leqslant i\leqslant l\leqslant M+N\ts)
$$ such that
$\la_{li}=\mu_i$ for every $l=1\lc M$. 
The space $V_\lm$ is the direct sum of the subspaces $V_\La$ in $V_\la$
where 
$\La$ runs through the set $\Slm$. We will assume that the set $\Slm$
is not empty. 

The action in $V_\la$ of the centralizer of the subalgebra $\UM$ in $\UMN$
preserves the subspace $V_\lm$ by the definition of this subspace. 
By Proposition 1.3 the image of the homomorphism 
$$
\pi_{M+N}\ts\ps^{\vp1}_\bl:\ts\YN\to\UMN;
\qquad
\bl=(M+1\lc M+N)
$$
commutes with the subalgebra $\UM$ in $\UMN$.
We will regard $V_\lm$as a module over the Yangian $\YN$ by making use of
this homomorphism.
Introduce the rational function 
$$
\al_{\mu}(u)=\prod_{i=1}^M\ts\frac
{\left(u+\mu_i-i+1\right)}{\left(u-i+1\right)}\ts;
$$
it can be also regarded as a formal power series in $u^{-1}$.
We have $\al_{M,\La}(u)=\al_{\mu}(u)$ for every $\La\in\Slm\ts$.

\proclaim{Proposition 2.2}
Suppose that $m\in\{1\lc N\}$ and $\La\in\Slm$.
The subspace $V_\La\subset V_\lm$ is an eigenspace for the coefficients
of the series
$A_m(u)$, the eigenvalues being the
respective coefficients of the series $\al_{M+m,\La}(u)/\al_{\mu}(u-m)$.
\endproclaim

\demo{Proof}
Let us apply Corollary 1.8 to $\bi=\bj=(1\lc m)$. Then 
we get the equality
$$
\ps^{\vp1}_\bl\bigl(A_m(u)\bigr)=A_{M+m}(u)/A_M(u-m);
\qquad
\bl=(M+1\lc M+N)
$$
in the algebra $\YMN[[u^{-1}]]$.
By applying to this equality the homomorphism $\pi_{M+N}$ and then
making use of Lemma 2.1 we obtain the required statement
\enddemos

\nt
The space $V_\lm$ inherits from $V_\la$ the $\ZZ$-grading.
Consider the $\ZZ$-grading $\operatorname{deg}$ on the algebra $\YN$
defined by \(1.999). 

\proclaim{Proposition 2.3}
The action of the algebra $\YN$ in the space $V_\lm$ is graded.
\endproclaim

\demo{Proof}
The action of the algebra $\UMN$ in the space $V_\la$ is a graded action.
The homomorphism $\pi_{M+N}:\YMN\to\UMN$ preserves the $\ZZ$-grading.
But the embedding $\ps^{\vp1}_\bl:\YN\to\YMN$ with 
$\bl=(M+1\lc M+N)$ also preserves the $\ZZ$-grading.
Indeed, we have 
$
\ps^{\vp1}_\bl=\si_{M+N}\ts\ph^{\vp1}_\bl\ts\si_M.
$
Here the embedding $\ph^{\vp1}_\bl:\YN\to\YMN$ preserves
$\operatorname{deg}$
by definition. But due to \(1.3) we have 
$\operatorname{deg}\cdot\ts\si_N=-\operatorname{deg}$
and similarly
$\operatorname{deg}\cdot\ts\si_{M+N}=-\operatorname{deg}$
\enddemos

\nt
Denote by $\La_\tz$ the array 
$(\ts\ka_{li}\ts|\ts 1\leqslant i\leqslant l\leqslant M+N\ts)$ where
$$
\kappa_{li}=
\cases
\mu_i&
\quad\text{if}\quad
l<M,
\\
\operatorname{min}\bigl(\ts\la_i,\mu_{i-l+M}\ts\bigr)&
\quad\text{if}\quad
l\geqslant M
\quad\text{and}\quad 
i>l-M,
\\
\la_i&
\quad\text{if}\quad
l>M
\quad\text{and}\quad 
i\leqslant l-M.
\endcases
\Tag{2.27}
$$

\proclaim{Lemma 2.4}
We have $\La_\tz\in\Slm\ts$.
\endproclaim

\demo{Proof}
For every $i\leqslant N$ we have $\ka_{M+N,i}=\la_i$ by definition.
Since the set $\Slm$ is non-empty, we have the inequality
$\la_i\leqslant\mu_{i-N}$ for
every $i>N$. Therefore $\ka_{M+N,i}=\la_i$ for $i>N$ also.
Similarly, for every $i\leqslant M$ we have $\la_i\geqslant\mu_i$ 
and $\ka_{Mi}=\mu_i$.
Since the sequences $\la$ and $\mu$ are non-increasing, we have the
inequalities
$\ka_{li}\geqslant\ka_{l-1,i}\geqslant\ka_{l,i-1}$ for all possible
$l$ and $i$
\enddemos 
 
\nt
Observe that for any scheme $\La\in\Slm$ we have the inequalities
$\la_{li}\leqslant\ka_{li}$ for all $l$ and $i$. Therefore the subspace
$V_{\La_\tz}$ in $V_\lm$ has the maximal degree, see \(2.07). Thus
we obtain the following corollary to Proposition 2.3.

\proclaim{Corollary 2.5}
The subspace $V_{\La_\tz}\subset V_\lm$ is annihilated by the action of any
element in $\YN$ of a positive degree.
\endproclaim

\nt
We will denote by $\Tlm$ the subset in $\Tn$ formed by the arrays $\La$
such that
$\la_{li}=\mu_i$ for each $l=1\lc M$ and $\la_{li}\leqslant\ka_{li}$ for
all $l$ and $i$.
We have $\al_{M,\La}(u)=\al_{\mu}(u)$ for every $\La\in\Tlm\ts$.

Let
$
\be=(\ts\be_1,\be_2,\dots)
$ 
and
$
\ga=(\ts\ga_1,\ga_2,\dots)
$
be any two non-increasing sequences of integers such that
$\be_i\geqslant\ga_i$ for each
$i=1,2,\dots$ and $\be_i=\ga_i$ for every  $i$
large  enough.
The {\it skew Young diagram} corresponding to the sequences
$\be,\ga$ is the set
$$
\{\ts(i,j)\in\ZZ^2\ts|\ts i\geqslant1,\ \be_i\geqslant j>\ga_i\ts\}
$$
We shall employ the usual graphic representation of a diagram:
a point $(i,j)\in\ZZ^2$
is represented by the unit box on the plane $\Bbb R^2$ with the centre
$(i,j)$,
the coordinates $i$ and $j$ on $\Bbb R^2$ increasing from top to bottom
and from left to
right respectively. The {\it content} of the box corresponding to $(i,j)$
is the
difference $j-i$.  Here is the diagram corresponding to the sequences
$$
\beta=(\ts6,\ts6,\ts4,\ts2,\ts0,\mi2,\mi3,\mi3,\ldots)
\quad\text{and}\quad
\gamma=(\ts3,\ts2,\ts2,\ts1,\mi3,\mi3,\ldots);
$$
in this diagram we have indicated
the content of the bottom box for every column.

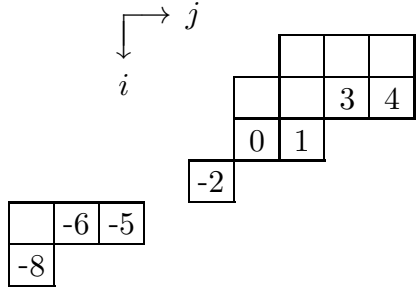
\vsk>
\vbox{
$$
\kern23.1pt\longrightarrow\,j\mv
$$
\vglue-26.2pt
$$
\vert\mv
$$
\vglue-28pt
$$
\bigr\downarrow\mv\kern0.1pt
$$
\vglue-16pt
$$
i\mv
$$
\vglue-40pt
$$
\phantom{\bx}
\phantom{\bx}
\phantom{\bx}
\phantom{\bx}
\phantom{\bx}
\phantom{\bx}
{\bx}
{\bx}
{\bx}
$$
\vglue-17.8pt
$$
\phantom{\bx}
\phantom{\bx}
\phantom{\bx}
\phantom{\bx}
\phantom{\bx}
{\bx}
{\bx}
{\bx}
{\bx}
$$
\vglue-17.8pt
$$
\phantom{\bx}
\phantom{\bx}
\phantom{\bx}
\phantom{\bx}
\phantom{\bx}
{\bx}
{\bx}
\phantom{\bx}
\phantom{\bx}
$$
\vglue-17.8pt
$$
\phantom{\bx}
\phantom{\bx}
\phantom{\bx}
\phantom{\bx}
{\bx}
\phantom{\bx}
\phantom{\bx}
\phantom{\bx}
\phantom{\bx}
$$
\vglue-17.8pt
$$
{\bx}
{\bx}
{\bx}
\phantom{\bx}
\phantom{\bx}
\phantom{\bx}
\phantom{\bx}
\phantom{\bx}
\phantom{\bx}
$$
\vglue-17.8pt
$$
{\bx}
\phantom{\bx}
\phantom{\bx}
\phantom{\bx}
\phantom{\bx}
\phantom{\bx}
\phantom{\bx}
\phantom{\bx}
\phantom{\bx}
$$
\vglue-83.5pt
$$
\kern100pt
3
\kern9pt
4
$$
\vglue-18pt
$$
\kern43pt
0
\kern9pt
1
$$
\vglue-17.8pt
$$
\mi2
$$
\vglue-17.8pt
$$
\phantom{\mi7}
\kern6pt
\mi6
\kern6pt
\mi5
\kern85pt
$$
\vglue-17.8pt
$$
\mi8
\kern113pt
$$
}
\vsk>
Now consider the scheme $\La_\tz\in\Slm\ts$. We have
$\ka_{M+m,M+m}\geqslant\ka_{M+N,M+N}=\la_{M+N}$ and
$\ka_{M+m,i}\geqslant\ka_{Mi}$ for all possible $m$ and $i$.
For every $m=0,1\lc N$  
denote by $\ka^{(m)}$ the skew Young diagram corresponding to the sequences
$$
\align
(\ts\ka_{M+m,1},\ka_{M+m,2}\ts&\lc\ka_{M+m,M+m},\la_{M+N},\la_{M+N}\dots\ts),
\\
(\ts\ka_{M1},\ka_{M2}\ts&\lc\ka_{MM},\la_{M+N},\la_{M+N}\dots\ts).
\endalign
$$
Thus the diagram $\ka^{(0)}$ is empty. The diagram $\ka^{(N)}$
will be denoted by $\la/\mu$.

\proclaim{Lemma 2.6}
Any column of the diagram $\ka^{(m)}$ consists of at most $m$ boxes.
For $m\geqslant1$
the diagram $\ka^{(m-1)}$ is obtained from $\ka^{(m)}$
by removing the bottom box from every column of the height $m$.
\endproclaim

\demo{Proof}
If $M=0$ then both $\ka^{(m)}$ and $\ka^{(m-1)}$ are usual Young diagrams.
The latter diagram is obtained from the former 
by removing the row $m$. Lemma 2.6 is
then evidently true. We will suppose that $M\geqslant1$.
If the top and the bottom boxes of a column in $\ka^{(m)}$ are $(i,j)$ and
$(i+l,j)$ respectively then $\ka_{Mi}<j\leqslant\ka_{M+m,i+l}$. But for
every 
$l\geqslant m$ we have the ineqialities
$\ka_{Mi}\geqslant\ka_{M+m,i+m}\geqslant\ka_{M+m,i+l}$ since
$\La_\tz\in\Slm\ts$.
This proves the first statement of Lemma 2.6.

To prove the second statement
observe that $\ka_{M+m-1,i}$ differs from $\ka_{M+m,i}$ only if
$i\geqslant m\geqslant 1$ and $\ka_{M+m,i}>\ka_{M,i-m+1}$. In this case
we have
$\ka_{M+m-1,i}=\ka_{M,i-m+1}$.
Each of the columns corresponding to $j$ with
$\ka_{M,i-m+1}<j\leqslant\ka_{M+m,i}$
has the height $m$. The box $(i,j)$ is then the bottom box of this column.
The row $i$ of $\ka^{(m-1)}$ is obtained from that of $\ka^{(m)}$
by removing each of these boxes
\enddemos

\nt
Note that the condition $\Slm\neq\varnothing$ on the sequences \(2.0744)
is equivalent to the condition that any column of the diagram $\la/\mu$
consists
of at most $N$ boxes. Now
for every $m=1\lc N-1$ introduce the polynomial in $u$
$$
P_{m,\la/\mu}(u)=\prod_k\ts(u+k)
\Tag{2.27333}
$$
where $k$ runs through the collection of the contents of the bottom boxes
in the columns of the height $m$ in the diagram $\la/\mu$ .
We will make use of the following observation.

\proclaim{Proposition 2.7}
For every $m=1\lc N-1$ we have the equality
$$
\frac
{
\al_{M+m+1,\La_\tz}(u)\ts\al_{M+m-1,\La_\tz}(u-1)
}{
\al_{M+m,\La_\tz}(u)\ts\al_{M+m,\La_\tz}(u-1)
}
=
\frac
{
P_{m,\la/\mu}(u-1)
}{
P_{m,\la/\mu}(u)
}\ts.
$$
\endproclaim

\demo{Proof}
For every $m=0\lc N$ introduce the polynomial
$$
Q_m(u)=\prod_k\ts(u+k)
$$
where $k$ runs through the collection of the contents of the bottom boxes
in the columns of the height $m$ in the diagram $\ka^{(m)}$. Then by
Lemma 2.6 we have
$$
P_{m,\la/\mu}(u)=Q_m(u)\ts\big/\ts Q_{m+1}(u+1);
\qquad
m=1\lc N-1.
\Tag{2.27111}
$$
On the other hand, by the definition \(2.0733) and again by Lemma 2.6
we have
$$
\frac
{
\al_{M+m,\La_\tz}(u)
}{
\al_{M+m-1,\La_\tz}(u)
}
=(u-M-m+1)\cdot
\frac
{
Q_m(u+1)
}{
Q_{m}(u)
}\ts;
\qquad
m=1\lc N.
\Tag{2.27222}
$$
By comparing \(2.27111) and \(2.27222) we obtain the required statement
\enddemos

\nt
The relations \(1.2) imply that for any $h\in\CC$ 
the assignment $T_{ij}(u)\mapsto T_{ij}(u+h)$ 
determines an automorphism of the algebra $\YN$;
here the series in $(u+h)^{-1}$ should be re-expanded in $u^{-1}$.
We will denote by 
$V_\lm(h)$ the $\YN$-module obtained from $V_\lm$ by the pullback through
this automorphism.

Let us now make use of the comultiplication \(1.135). For every $s=1\lc n$
fix some $h^{(s)}\in\CC$ along with a pair of non-increasing sequences
of integers
$$
\la^{(s)}=
\bigl(\ts\la_1^{(s)}\lc\la_{M^{(s)}}^{(s)},
\la_{{M^{(s)}}+1}^{(s)}\lc\la_{{M^{(s)}}+N}^{(s)}\ts\bigr)
\quad\text{and}\quad
\mu^{(s)}=\bigl(\ts\mu_1^{(s)}\lc\mu_{M^{(s)}}^{(s)}\ts\bigr)
$$
where ${M^{(s)}}\in\{\ts0,1,2,\ldots\ts\}$. Suppose that each of the sets
$S_{\la^{(s)},\mu^{(s)}}$ is not empty.
Consider the $\YN$-module
$$
W=V_{\la^{(1)},\mu^{(1)}}\bigl(h^{(1)}\bigr)
\ot\ldots\ot 
V_{\la^{(n)},\mu^{(n)}}\bigl(h^{(n)}\bigr).
\Tag{2.29}
$$
In Section 3 we will prove that the $\YN$-module $W$ is irreducible
provided $h^{(r)}-h^{(s)}\notin\ZZ$ for all $r\neq s$.
Theorem 2.9 below is contained in [C], Theorem 2.6.
However, we will give here the proof.  

For each $m=1\lc N$ and every collection of schemes
$\La^{(s)}\in\S_{\la^{(s)},\mu^{(s)}}$
where
$s=1\lc n$
introduce the rational function
$$
\chi_{\ts m,\La^{(1)}\lc\La^{(n)}}(u)=
\prod_{s=1}^n\ts\frac
{\al_{{M^{(s)}}+m,\La^{(s)}}\bigl(u+h^{(s)}\bigr)}
{\al_{\mu^{(s)}}\bigl(u-m+h^{(s)}\bigr)}\ts;
\Tag{2.2222}
$$
it can be also regarded as a formal power series in $u^{-1}$.

\proclaim{Proposition 2.8}
\!For any collection of schemes
$\!\La^{(s)}\in\S_{\la^{(s)},\mu^{(s)}}$
with $s=1\lc n$
there exists a non-zero vector $\xi$ in $W$ such that for every
$m=1\lc N$ we have
the equality
$$
A_m(u)\cdot\xi=\xi\cdot\chi_{\ts m,\La^{(1)}\lc\La^{(n)}}(u).
$$
\endproclaim

\demo{Proof}
We will derive this statement from Corollary 1.13.
Each of the vector spaces
$V_{\la^{(1)},\mu^{(1)}}\lc V_{\la^{(n)},\mu^{(n)}}$
is $\ZZ$-graded. Their tensor product acquires grading by the group $\ZZ^n$;
we have fixed the lexicographical ordering on the group $\ZZ^n$.

For each $s=1,\dots,\nomathbreak n$ and every
$\La^{(s)}\in\S_{\la^{(s)},\mu^{(s)}}$
choose a non-zero vector $\xi_{\La^{(s)}}$ in the subspace
$V_{\La^{(s)}}\subset V_{\la^{(s)},\mu^{(s)}}$. The vectors
$\xi_{\La^{(1)}}\ot\ldots\ot\xi_{\La^{(n)}}$
form a homogeneous basis in $W$.
By Corollary 1.13 and Proposition 2.3 
for every $m=1\lc N$
we have the equality of formal series in $u^{-1}$
$$
A_m(u)\cdot
\xi_{\La^{(1)}}\ot\ldots\ot\xi_{\La^{(n)}}=
\xi_{\La^{(1)}}\ot\ldots\ot\xi_{\La^{(n)}}\cdot
\chi_{\ts m,\La^{(1)}\lc\La^{(n)}}(u)\ts+\ts\Xi(u)
$$
where $\Xi(u)$ involves only the basic vectors in $W$ 
with the degrees smaller than that of
$\xi_{\La^{(1)}}\ot\ldots\ot\xi_{\La^{(n)}}$.
Proposition 2.8 now follows
\enddemos

\proclaim{Theorem 2.9}
Suppose that $h^{(r)}-h^{(s)}\notin\ZZ$ for all $r\neq s$. Then there is
a basis 
$$
\bigl(\ts\xi_{\La^{(1)}\lc\La^{(n)}}
\ts|\ts
\La^{(s)}\in\S_{\la^{(s)},\mu^{(s)}}\ts;\ s=1\lc n\ts)
$$
in the space $W$ such that for every $m=1\lc N$ we have the equality
$$
A_m(u)\cdot\xi_{\La^{(1)}\lc\La^{(n)}}=
\xi_{\La^{(1)}\lc\La^{(n)}}\cdot
\chi_{\ts m,\La^{(1)}\lc\La^{(n)}}(u).
\Tag{2.3}
$$                                                                                                                                                                                                                                                                                          
The equalities \(2.3) determine the vector
$\xi_{\La^{(1)}\lc\La^{(n)}}$ 
uniquely up to a scalar factor.
\endproclaim

\demo{Proof}
Since $h^{(r)}-h^{(s)}\notin\ZZ$ for all $r\neq s$,
the collection of the schemes
$$
\La^{(s)}\in\S_{\la^{(s)},\mu^{(s)}};
\quad
s=1\lc n
\nopagebreak
$$
can be uniquely restored from the the collection of the rational functions
$$
\chi_{\ts m,\La^{(1)}\lc\La^{(n)}}(u);
\quad
m=1\lc N.
$$
Thus the action in $W$ of the commutative algebra $\AN$  has a simple
spectrum
by Proposition 2.8.
The vector $\xi_{\La^{(1)}\lc\La^{(n)}}\in W$ can be determined as an
eigenvector
of this action satisfying the equality \(2.3) for every $m=1\lc N$
\enddemos

\nt
Let $V$ be an irreducible finite-dimensional module of the algebra
$\YN$. A basis in the vector space $V$ consisting of the eigenvectors of the
subalgebra $\AN$ in $\YN$ will be called a {\it Gelfand-Zelin basis}.
In Section 4
we will demonstrate that any module $V$ which admits such a basis,
can be obtained from
some module $W$ by applying to the algebra $\YN$ an automorphism of the form 
$\of$.
 
A non-zero vector $\xi\in V$ is called {\it singular} if it is annihilated
by all the
coef\-fi\-cients of the series $C_1(u)\lc C_{N-1}(u)$.
The vector $\xi$ is then unique
up to a scalar multiplier and is an eigenvector for the coefficients of the
series $A_1(u)\lc A_N(u)$; see  [D2], Theorem 2. Moreover, then
$$
\frac
{A_{m+1}(u)\ts A_{m-1}(u-1)}
{A_m(u)\ts A_m(u-1)}\cdot\xi
=
\frac{P_m(u-1)}{P_m(u)}\cdot\xi\ts;
\qquad
m=1\lc N-1
\Tag{2.31}
$$
for certain monic polynomials $P_1(u)\lc P_{N-1}(u)$ with the coefficients
in $\CC$.
These $N-1$ polynomials are called the {\it Drinfeld polynomials} of the
module $V$.
Every collection of $N-1$ monic polynomials arises in this way. The modules
with the
same Drinfeld polynomials may differ only by an automorphism of the algebra
$\YN$ of
the form $\of$. The following properties [D2] of the vector $\xi$ will be
used
in Section~4.

\proclaim{Proposition 2.10}
The vector $\xi$ is annihilated by every coefficient of the series
$T_{ij}(u)$ with
$i>j$. Furthermore, we have
$$
\frac{T_{m+1,m+1}(u-m)}{T_{m,m}(u-m)}\cdot\xi
=
\frac{P_m(u-1)}{P_m(u)}\cdot\xi\ts;
\qquad
m=1\lc N-1.
$$
\endproclaim

\nt
Denote by $\overline{V}$
the irreducible $\YN$-module obtained from $V$ by the pullback
through the automorphism $\tau_N$.
Let $\overline{P}_1(u)\lc\overline{P}_{N-1}(u)$ be the Drinfeld polynomials
of the module $\overline{V}$.
The next proposition will be also used in Section~4.

\proclaim{Proposition 2.11}
For each $m=1\lc N-1$ we have $\overline{P}_m(u)=P_{N-m}(u-m)$.
\endproclaim

\demo{Proof}
By the equalities 
\(1.655555) the automorphism $\tau_N$ preserves the subalgebra in $\YN$ 
generated by the coefficients of the series $C_1(u)\lc C_{N-1}(u)$ and
$A_N(u)$.
Therefore $V$ and $\overline{V}$ have the same singular vectors. The
required statement
now follows from the equalities \(1.666666) and \(2.31).
\enddemos

\nt
In the remainder of this section we will point out a singular vector in the
module $W$
and evaluate the corresponding Drinfeld polynomials provided $W$ is
irreducible. 
For every $s=1\lc n$ determine the scheme
$$
\La_\tz^{(s)}=
(\ts\ka_{li}^{(s)}\ts|\ts 1\leqslant i\leqslant l\leqslant M^{(s)}+N\ts)
\in\S_{\la^{(s)},\mu^{(s)}}
$$
in the way analogous to \(2.27): we put
$$
\kappa_{li}^{(s)}=
\cases
\mu_i^{(s)}&
\quad\text{if}\quad
l<M^{(s)},
\\
\operatorname{min}\bigl(\ts\la_i^{(s)},\mu_{i-l+M^{(s)}}^{(s)}\ts\bigr)&
\quad\text{if}\quad
l\geqslant M^{(s)}
\quad\text{and}\quad 
i>l-M^{(s)},
\\
\la_i^{(s)}&
\quad\text{if}\quad
l>M^{(s)}
\quad\text{and}\quad 
i\leqslant l-M^{(s)}.
\endcases
$$

\nt
Choose a non-zero vector $\xi_{\La_\tz^{(s)}}$ in the subspace
$V_{\La_\tz^{(s)}}\subset V_{\la^{(s)},\mu^{(s)}}$. 
Consider the vector
$\xi_0=\xi_{\La_\tz^{(1)}}\ot\ldots\ot\xi_{\La_\tz^{(n)}}$
in the space $W$.

\proclaim{Proposition 2.12}
The vector $\xi_0$ in $W$ is annihilated by all the coefficients of
the series 
$C_1(u)\lc C_{N-1}(u)$. For each $m=1\lc N$
we have the equality 
$$
A_m(u)\cdot
\xi_0=
\xi_0\ts\cdot\ts
\chi_{\ts m,\La_\tz^{(1)}\lc\La_\tz^{(n)}}(u).
\Tag{2.32}
$$
\endproclaim

\demo{Proof}
Denote by $\Xn$ the vector subspace in $\YN^{\ot n}[[u^{-1}]]$
spanned by all homogeneous elements whose degree in $\ZZ^n$ has at least
one
positive component. Since $\deg C_m(u)=1$, by the definition of the
comultiplication
\(1.135) we have
$$
\De^{(n)}\bigl(C_m(u)\bigr)\in\Xn\ts;\qquad m=1\lc N-1.
$$
Therefore $C_m(u)\cdot\xi_0=0$ by Corollary 2.5. Furthermore,
by Corollary 1.13 we have
$$
\De^{(n)}\bigl(A_m(u)\bigr)\in\bigl(A_m(u)\bigr)^{\ot n}+\Xn\ts;
\qquad m=1\lc N-1.
$$
Thus by Proposition 2.2 the vector 
$\xi_0$ satisfies \(2.32) for every $m=1\lc N$
\enddemos

\nt
We now obtain the following corollary to Propositions 2.7 and 2.12.

\proclaim{Corollary 2.13}
Suppose that the module $W$ is irreducible. 
Then the vector $\xi_0$ in $W$ is singular. 
The Drinfeld polynomials of the module $W$ are
$$
\prod_{s=1}^n\ts P_{m,\la^{(s)}/\mu^{(s)}}\bigl(u+h^{(s)}\bigr)\ts;
\qquad
m=1\lc N-1.
\Tag{2.9999}
$$
\endproclaim

\nt
This fact will be substantially used in Section 4.

\bigskip\bigskip
\Section{Action of the Drinfeld generators} 

\nt
In this section we will keep fixed the sequences
$\la^{(1)},\mu^{(1)}\lc\la^{(n)},\mu^{(n)}$ along with the  parameters
$h^{(1)}\lc h^{(n)}\in\CC$. From now on until the end of this section we
will be assuming
that $h^{(r)}-h^{(s)}\notin\ZZ$ for all $r\neq s$.
For each $m=0\lc N$ denote
$$
\align
\rho_m(u)
&\ts=\ts
\prod_{s=1}^n
\left(\ts
\al_{\mu^{(s)}}\bigl(u-m+h^{(s)}\bigr)
\prod_{i=1}^{M^{(s)}+m}\ts\bigl(u-i+1+h^{(s)}\bigr)
\right)
\\
&\ts=\ts
\prod_{s=1}^n
\left(\ts
\prod_{i=1}^m\ts\bigl(u-i+1+h^{(s)}\bigr)
\prod_{i=1+m}^{M^{(s)}+m}\ts\bigl(u+\mu_{i-m}^{(s)}-i+1+h^{(s)}\bigr)
\right);
\endalign
$$
it is a polynomial in $u$ of the degree
$$
r_m=M^{(1)}+\ldots+M^{(n)}+m\ts n.
$$
Note that for each $m=1\lc N-1$ we have the equality
$$
\rho_m(u)\ts\rho_m(u-1)=\rho_{m+1}(u)\ts\rho_{m-1}(u-1).
\Tag{3.0}
$$
Consider the action of the coefficients of the formal Laurent series in
$u^{-1}$
$$
\rho_m(u)\ts A_m(u)
\quad\text{and}\quad
\rho_m(u)\ts B_m(u),
\quad
\rho_m(u)\ts C_m(u),
\quad
\rho_m(u)\ts D_m(u)
\Tag{3.1}
$$
in the $\YN$-module $W$.
Denote by
$a_m(u)$
and
$b_m(u),c_m(u),d_m(u)$
the series with the coefficients in $\End(W)$ corresponding to the series
\(3.1).
Observe that $a_m(u)$ is a polynomial in $u$ by Theorem 2.9.

\proclaim{Proposition 3.1}
For every $m=1\lc N-1$ these $b_m(u),c_m(u)$ and $d_m(u)$ are polynomial
in $u$.
The degrees of $b_m(u)$ and $c_m(u)$ do not exceed $r_m-1$.
\endproclaim

\demo{Proof}
Due to Proposition 1.11 each of the images
$$
\De^{(n)}\bigl(B_m(u)\bigr),\ts
\De^{(n)}\bigl(C_m(u)\bigr),\ts
\De^{(n)}\bigl(D_m(u)\bigr)
\in\YN[[u^{-1}]]
$$ 
is a finite sum of the series of the form 
$$
Q_{\bk^{(0)}\bk^{(1)}}(u)\ot
Q_{\bk^{(1)}\bk^{(2)}}(u)\ot\ldots\ot
Q_{\bk^{(n-1)}\bk^{(n)}}(u)
$$
where $\bk^{(0)},\bk^{(1)}\lc\bk^{(n)}$ are increasing sequences of the
integers
$1\lc N$ of the length $m$.
For each $\bk^{(s)}=\bigl(k_1^{(s)}\lc k_m^{(s)}\bigr)$ denote
$$
\bold{\bar k}^{(s)}=
\bigl(
1\lc M^{(s)},
M^{(s)}+k_1^{(s)}\lc M^{(s)}+k_m^{(s)}
\bigr).
$$
Due to Corollary 1.8 for each $s=1\lc n$
the action of the coefficients of the series $Q_{\bk^{(s-1)}\bk^{(s)}}(u)$
in
the module $V_{\la^{(s)},\mu^{(s)}}\bigl(h^{(s)}\bigr)$
can be defined as that of the coefficients
\nolinebreak
of
$$
\pi_{M^{(s)}+N}
\bigl(
Q_{\bold{\bar k}^{(s-1)}\bold{\bar k}^{(s)}}
\bigl(u+h^{(s)}\bigr)
\cdot
A_M^{(s)}
\bigl(u-m+h^{(s)}\bigr)^{-1}
\bigr)
\in
\operatorname{U}(\gl_{M^{(s)}+N})[[u^{-1}]]
$$
in the subspace $V_{\la^{(s)},\mu^{(s)}}\subset V_{\la^{(s)}}$. But the
product
$$
\pi_{M^{(s)}+N}
\bigl(
Q_{\bold{\bar k}^{(s-1)}\bold{\bar k}^{(s)}}
\bigl(u+h^{(s)}\bigr)
\bigr)
\cdot
\prod_{i=1}^{M^{(s)}+m}\ts\bigl(u-i+1+h^{(s)}\bigr)
$$
is a polynomial in $u$ by the definition of the homomorphism
$\pi_{M^{(s)}+N}$. 
By Lemma 2.1 the coefficients of 
$$
\pi_{M^{(s)}+N}\bigl(A_M^{(s)}\bigl(u-m+h^{(s)}\bigr)\bigr)
\in
\operatorname{U}(\gl_{M^{(s)}+N})[[u^{-1}]]
$$
act in $V_{\la^{(s)},\mu^{(s)}}$ as the respective coefficients of the
series
$\al_{\mu^{(s)}}\bigl(u-m+h^{(s)}\bigr)$.
Therefore $b_m(u),c_m(u)$ and $d_m(u)$ are polynomial in $u$.
Now observe that
$$
B_m(u),C_m(u)\in\YN[[u^{-1}]]\cdot u^{-1}
$$
by the definition \(1.605).
Since the degree of the polynomial $\rho_m(u)$ equals $r_m$
the degrees of the polynomials $b_m(u)$ and $c_m(u)$ do not exceed $r_m-1$
\enddemos

\nt
In particular, $a_m(u)$ and $b_m(u),c_m(u),d_m(u)$ can be evaluated
at any point
$u\in\CC$. Let us now fix for every $s=1\lc n$ a scheme
$$
\La^{(s)}=
(\ts\la_{li}^{(s)}\ts|\ts 1\leqslant i\leqslant l\leqslant M^{(s)}+N\ts)
\in\T_{\la^{(s)},\mu^{(s)}}\ts.
$$
For every $m=0\lc N$ consider the polynomial
$$
\align
\varpi_{m,\La^{(1)}\lc\La^{(n)}}(u)
&=
\rho_m(u)\ts\ts\cdot\ts
\chi_{\ts m,\La^{(1)}\lc\La^{(n)}}(u)
\\
&=\ts
\prod_{s=1}^n
\prod_{i=1}^{M^{(s)}+m}\ts\bigl(u+\la_{M^{(s)}+m,i}^{(s)}-i+1+h^{(s)}\bigr);
\endalign
$$
see the definition \(2.2222). Note that all the $r_m$ zeroes of 
this polynomial
$$
\nu_{mi}^{(s)}=i-\la_{M^{(s)}+m,i}^{(s)}-1-h^{(s)}
$$
are pairwise distinct due to our assumption on the parameters
$h^{(1)}\lc h^{(n)}$.
Also that we have the equality
$
\varpi_{0,\La^{(1)}\lc\La^{(n)}}(u)=\rho_0(u).
$

Endow the set of the pairs $(l,j)$ where $l=1\lc N$ and $j=1,2,\dots$
with the following relation of precedence:
$(m,i)\prec (l,j)$ if $i<j$ or $i=j$ and $m>l$. 
Let $\xi_0$ be the singular vector in $W$ constructed in the end of
Section 2. 
Consider the vector in $W$
$$
\xi=\ts\prod_{(l,j)}^\rightarrow\ 
\left(\ts
\prod_{(r,t)}\ts
b_l(\nu_{lj}^{(r)}-t)
\right)\cdot\xi_0
\Tag{3.2}
$$
where the product in brackets is taken over the set $\Cal F$ of all pairs 
$(r,t)$ such that
$$
r=1\lc n
\quad\text{and}\quad
j\leqslant M^{(r)}+l\ts;
\qquad
t=1\lc \ka_{M^{(r)}+l,j}^{(r)}-\la_{M^{(r)}+l,j}^{(r)}.
$$
Here for each fixed $l$ the elements
$b_l(\nu_{lj}^{(r)}-t)\in\End(W)$ commute because of the
relation \(1.9). The products in brackets do not commute with each other
in general. We
arrange them from the left to the right according to the above relation of
precedence for the pairs $(l,j)$.

\proclaim{Theorem 3.2}
For every $m=1\lc N$ we have the equality
$$
a_m(u)\cdot\xi=\varpi_{m,\La^{(1)}\lc\La^{(n)}}(u)\cdot\xi.
$$
\endproclaim 

\demo{Proof}
We will employ the induction on the number of the factors
$b_l(\nu_{lj}^{(r)}-t)$
in \(3.2).
If there is no factors then $\La^{(r)}=\La^{(r)}_0$ for every $r=1\lc n$
and $\xi=\xi_0$. In
particular, the required equality then holds by Proposition 2.12.

Assume that the product in \(3.2) contains at least one factor.
Let $(l,j)$ be the minimal pair such that the corresponding set $\Cal F$
is not empty. Let us fix an arbitrary $r\in\{1\lc n\}$ such that
$$
\la_{M^{(r)}+l,j}^{(r)}<\ka_{M^{(r)}+l,j}^{(r)}\ts.
$$
Denote by
$\Om^{(r)}$ the array obtained from $\La^{(r)}$ by increasing
the $(M^{(r)}+l,j)$-entry by $1\ts$. Then
$\Om^{(r)}\in\T_{\la^{(r)},\mu^{(r)}}$ and
$$
\xi=b_l(\nu_{lj}^{(r)}-1)\cdot\eta
\Tag{3.3}
$$
where
the vector $\eta$ is determined in the way analogous to \(3.2)  
by the sequence of schemes 
$\La^{(1)}\lc\Om^{(r)}\lc
\La^{(n)}$ instead of $\La^{(1)}\lc\La^{(r)}\lc\La^{(n)}$.
If $l\neq m$ then
$$
\varpi_{m,\La^{(1)}\lc\Om^{(r)}\lc\La^{(n)}}(u)=
\varpi_{m,\La^{(1)}\lc\La^{(r)}\lc\La^{(n)}}(u).
$$
By the relation \(1.7) and by the inductive assumption we get
$$
\gather
a_m(u)\cdot\xi=
a_m(u)\ts b_l(\nu_{lj}^{(r)}-1)\cdot\eta=
b_l(\nu_{lj}^{(r)}-1)\ts a_m(u)\cdot\eta=
\\
b_l(\nu_{lj}^{(r)}-1)\ts\varpi_{m,\La^{(1)}\lc\Om^{(r)}\lc\La^{(n)}}(u)
\cdot\eta=
\varpi_{m,\La^{(1)}\lc\La^{(r)}\lc\La^{(n)}}(u)\cdot\xi.
\endgather
$$

Now suppose that $l=m\ts$; then by the definition of $\Om^{(r)}$ we have
$$
\varpi_{m,\La^{(1)}\lc\Om^{(r)}\lc\La^{(n)}}(u)=
\frac{u-\nu_{mj}^{(r)}+1}{u-\nu_{mj}^{(r)}}\ts\ts
\varpi_{m,\La^{(1)}\lc\La^{(r)}\lc\La^{(n)}}(u).
$$
In particular, by the inductive assumption we then have
$$
a_m(\nu_{mj}^{(r)}-1)\cdot\eta=
\varpi_{m,\La^{(1)}\lc\Om^{(r)}\lc\La^{(n)}}(\nu_{mj}^{(r)}-1)\cdot\eta=0.
$$
Therefore by the relation \(1.10) and again by the inductive assumption
we get
$$
\align
a_m(u)\cdot\xi&=
a_m(u)\ts b_m(\nu_{mj}^{(r)}-1)\cdot\eta
\\
&=\frac{u-\nu_{mj}^{(r)}}{u-\nu_{mj}^{(r)}+1}\ts\ts
b_m(\nu_{mj}^{(r)}-1)\ts a_m(u)\cdot\eta
\\
&=\frac{u-\nu_{mj}^{(r)}}{u-\nu_{mj}^{(r)}+1}\ts\ts
b_m(\nu_{mj}^{(r)}-1)\ts\varpi_{m,\La^{(1)}\lc\Om^{(r)}\lc\La^{(n)}}(u)
\cdot\eta
\\
&=
\varpi_{m,\La^{(1)}\lc\La^{(r)}\lc\La^{(n)}}(u)\cdot\xi\ts\quad\square
\endalign
$$
\enddemo

\nt
Therefore the vector $\xi$ determined by \(3.2) satisfies the equality
\(2.3) for every
$m=1\lc N$. By Theorem 2.9 we can choose $\xi_{\La^{(1)}\lc\La^{(n)}}=\xi$
as soon as we
prove that $\xi\neq0$. We will describe first the action on the vector
$\xi$
of the coefficients of the polynomials $b_m(u)$ and $c_m(u)$ for
$m=1\lc N-1$.
Then we will prove that the vector $\xi$ does not vanish. 

\proclaim{Proposition 3.3}
If $\La^{(r)}\notin\S_{\la^{(r)},\mu^{(r)}}$ for some $r\in\{1\lc n\}$
then $\xi=0$.
\endproclaim 

\demo{Proof}
As well as in the proof of Theorem 3.2 we will employ the induction on the
number of the factors $b_l(\nu_{lj}^{(r)}-t)$ in \(3.2).
If there is no factors then $\La^{(r)}=\La^{(r)}_0$ for every $r=1\lc n$
and we have nothing to prove.

Assume that the product in \(3.2) contains at least one factor.
Let $(l,j)$ be the minimal pair such that the corresponding set $\Cal F$
is not empty. Fix an arbitrary $r\in\{1\lc n\}$ such that
$$
\la_{M^{(r)}+l,j}^{(r)}<\ka_{M^{(r)}+l,j}^{(r)}\ts.
$$
Introduce the array $\Om^{(r)}\in\T_{\la^{(r)},\mu^{(r)}}$ and determine
the
vector $\eta\in W$ in the same way as it was done in the proof of
Theorem 3.2.
Then we have the
equality \(3.3). Suppose that $\Om^{(r)}\in\S_{\la^{(r)},\mu^{(r)}}$
and $\La^{(s)}\in\S_{\la^{(s)},\mu^{(s)}}$ for any $s\neq r$.
Otherwise $\eta=0$ by the inductive assumption so that $\xi=0$ due
to \(3.3).

By Theorem 3.2 the vector $\xi$ is an
eigenvector for the coefficients of the polynomials
$a_1(u),\ldots,a_{n-1}(u)$.
Therefore by Theorem 2.9
$$
\varpi_{m,\La^{(1)}\lc\La^{(r)}\lc\La^{(n)}}(u)=
\varpi_{m,\Up^{(1)}\lc\Up^{(r)}\lc\Up^{(n)}}(u);
\qquad
m=1\lc N
$$
for certain schemes
$$
\Up^{(r)}\in\S_{\la^{(r)},\mu^{(r)}}\ts;
\qquad
r=1\lc n.
$$
Consider the roots
$\nu_{mi}^{(r)}$
of the polynomial $\varpi_{m,\La^{(1)}\lc\La^{(r)}\lc\La^{(n)}}(u)$.
Since the array
$\Om^{(r)}\!\in\nomathbreak\S_{\la^{(r)},\mu^{(r)}}$ we have the
inequalities
$$
\gather
\nu_{m1}^{(r)}+h^{(r)}<\nu_{m2}^{(r)}+h^{(r)}<
\ldots<\nu_{m,M^{(r)}+m}^{(r)}+h^{(r)}
\qquad\text{if}\quad
m\neq l;
\\
\nu_{l1}^{(r)}+h^{(r)}<\ldots<\nu_{lj}^{(r)}+h^{(r)}
\leqslant\nu_{l,j+1}^{(r)}+h^{(r)}
<\ldots<\nu_{l,M^{(r)}+l}^{(r)}+h^{(r)}\ts.
\endgather
$$
Moreover, we have $h^{(s)}-h^{(r)}\notin\ZZ$ for any $s\neq r$ by our
assumption.
Therefore the array $\La^{(r)}$ can be uniquely restored from the 
collection
of the
polynomials
$$
\varpi_{m,\La^{(1)}\lc\La^{(r)}\lc\La^{(n)}}(u);
\qquad
m=1\lc N.
$$
Thus $\La^{(r)}=\Up^{(r)}\in\S_{\la^{(r)},\mu^{(r)}}$
and the Proposition 3.3 is proved
\enddemos

\nt
  From now on we will suppose that $\La^{(s)}\in\S_{\la^{(s)},\mu^{(s)}}$
for each $s=1\lc n$.  Let an index $m\in\{1\lc N-1\}$ be fixed.
By Proposition 3.1 to determine the polynomials 
$b_m(u)\cdot\xi$ and $c_m(u)\cdot\xi$ it suffices to evaluate them at
$r_m$ points.
We will evaluate these polynomials at $u=\nu_{mi}^{(s)}$ where 
$s=1\lc n$ and $i=1\lc {M^{(s)}+m}\ts$.

Let the indices $s$ and $i$ be fixed. 
Denote by $\La_\tp^{(s)}$ the array obtained from $\La^{(s)}$ by increasing
the $(M^{(s)}+m,i)$-entry by $1\ts$.
If $\La_\tp^{(s)}\in\S_{\la^{(s)},\mu^{(s)}}$
then denote by $\xi_+$ the vector in $W$ determined in the way
analogous to \(3.2)  
by the sequence of schemes 
$\La^{(1)}\lc\La_\tp^{(s)}\lc\La^{(n)}$ instead of
$\La^{(1)}\lc\La^{(s)}\lc\La^{(n)}$.
Denote by $\ga_{mi}^{(s)}$ the product
$$
\ \
\prod_{r=1}^n\ts
\prod_{j=1}^{M^{(r)}+m+1}
\ts
\cases
\bigl(
\ka_{M^{(r)}+m+1,j}^{(r)}-\la_{M^{(s)}+m,i}^{(s)}+
i-j+
h^{(r)}-h^{(s)}
\bigr)
&\quad\text{if}\quad
j\leqslant i
\\
\bigl(
\la_{M^{(r)}+m+1,j}^{(r)}-\la_{M^{(s)}+m,i}^{(s)}+
i-j+
h^{(r)}-h^{(s)}
\bigr)
&\quad\text{if}\quad
j>i\ts,
\endcases
$$
multiplied by the product
$$
\ 
\prod_{r=1}^n\ts
\prod_{j=1}^{M^{(r)}+m-1}
\ts
\cases
\bigl(
\ka_{M^{(r)}+m-1,j}^{(r)}-\la_{M^{(s)}+m,i}^{(s)}+
i-j-1+
h^{(r)}-h^{(s)}
\bigr)
&\ \text{if}\quad
j<i
\\
\bigl(
\la_{M^{(r)}+m-1,j}^{(r)}-\la_{M^{(s)}+m,i}^{(s)}+
i-j-1+
h^{(r)}-h^{(s)}
\bigr)
&\ \text{if}\quad
j\geqslant i\ts.
\endcases
$$

\proclaim{Lemma 3.4}
If $\La_\tp^{(s)}\in\S_{\la^{(s)},\mu^{(s)}}$ then $\ga_{mi}^{(s)}\neq0$.
\endproclaim

\demo{Proof}
Any factor in the product $\ga_{mi}^{(s)}$ corresponding to $r\neq s$ 
is not zero since $h^{(r)}-h^{(s)}\notin\ZZ$ then. Since
$\La^{(s)},\La_\tz^{(s)},\La_\tp^{(s)}\in\S_{\la^{(s)},\mu^{(s)}}$
we have the inequalities
$$
\alignat2
&
\ka_{M^{(s)}+m+1,j}^{(s)}
\geqslant
\ka_{M^{(s)}+m,j}^{(s)}
\geqslant
\ka_{M^{(s)}+m,i}^{(s)}
>
\la_{M^{(s)}+m,i}^{(s)}
&&
\qquad\text{if}\quad j\leqslant i\ts,
\\
&
\la_{M^{(s)}+m+1,j}^{(s)}
\leqslant
\la_{M^{(s)}+m,j-1}^{(s)}
\leqslant
\la_{M^{(s)}+m,i}^{(s)}
&&
\qquad\text{if}\quad j> i\ts,
\\
&
\ka_{M^{(s)}+m-1,j}^{(s)}
\geqslant
\ka_{M^{(s)}+m,j+1}^{(s)}
\geqslant
\ka_{M^{(s)}+m,i}^{(s)}
>
\la_{M^{(s)}+m,i}^{(s)}
&&
\qquad\text{if}\quad j< i\ts,
\\
&
\la_{M^{(s)}+m-1,j}^{(s)}
\leqslant
\la_{M^{(s)}+m,j}^{(s)}
\leqslant
\la_{M^{(s)}+m,i}^{(s)}
&&
\qquad\text{if}\quad j\geqslant i\ts.
\endalignat
$$
These inequalities show that
any factor in the product $\ga_{mi}^{(s)}$ corresponding to $r=s$ 
is also non-zero
\enddemos

\proclaim{Theorem 3.5}
We have 
$$
c_m\bigl(\nu_{mi}^{(s)}\bigr)\cdot\xi=
\cases
-\ts\ts\ga_{mi}^{(s)}\cdot\xi_+\quad
&\text{if}\quad\La_\tp^{(s)}\in\S_{\la^{(s)},\mu^{(s)}}\ts;
\\
\qquad0\quad
&\text{otherwise}.
\endcases
$$
\endproclaim

\demo{Proof}
If $\La^{(r)}=\La^{(r)}_0$ for every $r=1\lc n$ then
$\La_\tp^{(s)}\notin\S_{\la^{(s)},\mu^{(s)}}$ in particular.
On the other hand then
$\xi=\xi_0$ and by Proposition 2.12 we have
$
c_m(\nu_{mi}^{(s)})\cdot\xi=0\ts.
$

Assume that the product in \(3.2) contains at least one factor.
Let $(l,j)$ be the minimal pair such that the corresponding set $\Cal F$
is not empty. Fix an arbitrary $r\in\{1\lc n\}$ such that
$$
\la_{M^{(r)}+l,j}^{(r)}<\ka_{M^{(r)}+l,j}^{(r)}\ts.
$$
Introduce the array $\Om^{(r)}\in\T_{\la^{(r)},\mu^{(r)}}$ and determine
the
vector $\eta\in W$ in the same way as it was done in the proof of
Theorem 3.2.
Then we have the
equality \(3.3). 
Moreover, then $\Om^{(r)}\in\S_{\la^{(r)},\mu^{(r)}}$
since $\La^{(r)}\in\S_{\la^{(r)},\mu^{(r)}}$ by our assumption.

Let us prove first that the following equality holds for any
$\ts(l,j,r)\neq(m,i,s)$:
$$ 
c_m\bigl(\nu_{mi}^{(s)}\bigr)\ts b_l\bigl(\nu_{lj}^{(r)}-1\bigr)
\cdot\eta\ts=
b_l\bigl(\nu_{lj}^{(r)}-1\bigr)\ts c_m\bigl(\nu_{mi}^{(s)}\bigr)
\cdot\eta.
\Tag{3.4}
$$
If $l\neq m$ then we obtain \(3.4) directly from the relation \(1.8).
But if $l=m$
then $\nu_{mj}^{(r)}-1\neq\nu_{mi}^{(s)}$. Indeed, for $r\neq s$
it follows from
the assumption $h^{(r)}-h^{(s)}\notin\ZZ$.
Since $\Om^{(r)}\in\S_{\la^{(r)},\mu^{(r)}}$
we also have $\nu_{mj}^{(r)}-1\neq\nu_{mi}^{(r)}$. 
Due to Theorem 3.2 we then also have the equalities 
$$
\align
a_m\bigl(\nu_{mj}{(r)}-1\bigr)&\cdot\eta=
\varpi_{m,\La^{(1)}\lc\Om^{(r)}\lc\La^{(n)}}\bigl(\nu_{mj}^{(r)}-1\bigr)
\cdot\eta=0,
\\
a_m\bigl(\nu_{mi}{(s)}\bigr)&\cdot\eta=
\varpi_{m,\La^{(1)}\lc\Om^{(r)}\lc\La^{(n)}}\bigl(\nu_{mi}^{(s)}\bigr)
\cdot\eta=0.
\endalign
$$
Therefore by the relation \(1.11) we again obtain that
$$
c_m\bigl(\nu_{mi}^{(s)}\bigr)\ts b_m\bigl(\nu_{mj}^{(r)}-1\bigr)
\cdot\eta=
b_m\bigl(\nu_{mj}^{(r)}-1\bigr)\ts c_m\bigl(\nu_{mi}^{(s)}\bigr)
\cdot\eta\ts.
$$

If $\la_{M^{(s)}+m,i}^{(s)}=\ka_{M^{(s)}+m,i}^{(s)}$ then
$\La_\tp^{(s)}\notin\S_{\la^{(s)},\mu^{(s)}}\ts$.
On the other hand, by applying the equality \(3.4) repeatedly we then get
$$
\align
c_m\bigl(\nu_{mi}^{(s)}\bigr)\cdot\xi
&=
c_m\bigl(\nu_{mi}^{(s)}\bigr)\ts
\ts\prod_{(l,j)}^\rightarrow\ 
\left(\ts
\prod_{(r,t)}\ts
b_l(\nu_{lj}^{(r)}-t)
\right)\cdot\xi_0
\\
&=
\ts\prod_{(l,j)}^\rightarrow\ 
\left(\ts
\prod_{(r,t)}\ts
b_l(\nu_{lj}^{(r)}-t)
\right)\ts\ts
c_m\bigl(\nu_{mi}^{(s)}\bigr)
\cdot\xi_0=0
\endalign
$$
as we have claimed. Now we will assume that
$\la_{M^{(s)}+m,i}^{(s)}<\ka_{M^{(s)}+m,i}^{(s)}\ts$.

For every $r=1\lc n$ consider the array  
$\Up^{(r)}$ the array obtained from $\La^{(r)}$ by changing 
all the $(M^{(r)}+l,j)$-entries corresponding to the pairs
$(l,j)\prec(m,i)$
for $\ka_{M^{(r)}+l,j}^{(r)}$ and also by increasing
the $(m,i)$-entry by $1$ if $r=s$.
Then
$\Up^{(r)}\in\T_{\la^{(r)},\mu^{(r)}}$ for any $r$. Determine the
vector $\zeta$ in the way analogous to \(3.2)  
by the sequence of schemes 
$\Up^{(1)}\lc\Up^{(n)}$ instead of $\La^{(1)}\lc\La^{(n)}$.
Due to Theorem 3.2 we then have
$$
a_m\bigl(\nu_{mi}^{(s)}-1\bigr)\cdot\zeta=
\varpi_{m,\Up^{(1)}\lc\Up^{(n)}}\bigl(\nu_{mi}^{(s)}-1\bigr)
\cdot\zeta=0\ts.
$$
We then also have
$
\xi=p\ts b_m\bigl(\nu_{mi}^{(s)}-1\bigr)\cdot\zeta
$
where
$$
p=\ts\prod_{(l,j)\prec(m,i)}^\rightarrow\ 
\left(\ts
\prod_{(r,t)}\ts
b_l(\nu_{lj}^{(r)}-t)
\right).
$$
Therefore by applying the equality \(3.4) repeatedly and using
Proposition 1.10
along with the equality \(3.0) we get
$$
\gather
c_m\bigl(\nu_{mi}^{(s)}\bigr)\cdot\xi
=
p\ts c_m\bigl(\nu_{mi}^{(s)}\bigr)\ts b_m\bigl(\nu_{mi}^{(s)}-1\bigr)
\cdot\zeta
=-\ts
p\ts a_{m+1}\bigl(\nu_{mi}^{(s)}\bigr)\ts
a_{m-1}\bigl(\nu_{mi}^{(s)}-1\bigr)\cdot\zeta=
\\
-\ts
\varpi_{m+1,\Up^{(1)}\lc\Up^{(n)}}\bigl(\nu_{mi}^{(s)}\bigr)
\ts
\varpi_{m-1,\Up^{(1)}\lc\Up^{(n)}}\bigl(\nu_{mi}^{(s)}-1\bigr)
\ts p\cdot\zeta
=-\ts\ga_{mi}^{(s)}\ts p\cdot\zeta\ts.
\endgather
$$
Here if $\La_\tp^{(s)}\in\S_{\la^{(s)},\mu^{(s)}}$ then
$p\cdot\zeta=\xi_+$ by
definition. If $\La_\tp^{(s)}\notin\S_{\la^{(s)},\mu^{(s)}}$ then
$p\cdot\zeta=0$
by Proposition 3.3
\enddemos

\proclaim{Proposition 3.6}
If $\La^{(r)}\in\S_{\la^{(r)},\mu^{(r)}}$ for every $r=1\lc n$ then
$\xi\neq0$.
\endproclaim 

\demo{Proof}
As well as in the proof of Theorem 3.2 we will employ the induction on the
number of the factors $b_l(\nu_{lj}^{(r)}-t)$ in \(3.2).
If there is no factors then $\La^{(r)}=\La^{(r)}_0$ for every $r=1\lc n$
and $\xi=\xi_0\neq0$.

Assume that the product in \(3.2) contains at least one factor.
Let $(l,j)$ be the minimal pair such that the corresponding set $\Cal F$
is not empty. Fix an arbitrary $r\in\{1\lc n\}$ such that
$$
\la_{M^{(r)}+l,j}^{(r)}<\ka_{M^{(r)}+l,j}^{(r)}\ts.
$$
Introduce the array $\Om^{(r)}\in\T_{\la^{(r)},\mu^{(r)}}$ and determine
the
vector $\eta\in W$ in the same way as it was done in the proof of
Theorem 3.2.
Then we have the equality \(3.3).
Since $\La^{(r)}\in\S_{\la^{(r)},\mu^{(r)}}$ we also have
$\Om^{(r)}\in\S_{\la^{(r)},\mu^{(r)}}$. Therefore $\eta\neq0$ by the
inductive
assumtion. On the other hand, by Theorem 3.5
we then have
$$
c_l\bigl(\nu_{lj}^{(r)}\bigr)\cdot\xi=\ga_{lj}^{(r)}\cdot\eta
$$
where $\ga_{lj}^{(r)}\neq0$ due to Lemma 3.4. Therefore $\xi\neq0$
\enddemos

\nt
Now consider the array $\La_\tm^{(s)}$ obtained from $\La^{(s)}$ by
decreasing
the $(M^{(s)}+m,i)$-entry by $1\ts$.
If $\La_\tm^{(s)}\in\S_{\la^{(s)},\mu^{(s)}}$
then denote by $\xi_-$ the vector in $W$ determined in the way analogous
to \(3.2)  
by the sequence of schemes 
$\La^{(1)}\lc\La_\tm^{(s)}\lc\La^{(n)}$ instead
of $\La^{(1)}\lc\La^{(s)}\lc\La^{(n)}$.
Denote by $\be_{mi}^{(s)}$ the product
$$
\align
&
\prod_{r=1}^n\ts\ts
\prod_{j=1}^{\operatorname{min}(i,M^{(r)}+m+1)}
\ts\ts
\frac
{
\la_{M^{(r)}+m+1,j}^{(r)}-\la_{M^{(s)}+m,i}^{(s)}+
i-j+1+
h^{(r)}-h^{(s)}
}{
\ka_{M^{(r)}+m+1,j}^{(r)}-\la_{M^{(s)}+m,i}^{(s)}+
i-j+1+
h^{(r)}-h^{(s)}
}
\\
\times\ts\ts
&
\prod_{r=1}^n\ts\ts
\prod_{j=1}^{\operatorname{min}(i,M^{(r)}+m)-1}
\ts\ts
\frac
{
\la_{M^{(r)}+m-1,j}^{(r)}-\la_{M^{(s)}+m,i}^{(s)}+
i-j+
h^{(r)}-h^{(s)}
}{
\ka_{M^{(r)}+m-1,j}^{(r)}-\la_{M^{(s)}+m,i}^{(s)}+
i-j+
h^{(r)}-h^{(s)}
}\ts.
\endalign
$$

\proclaim{Lemma 3.7}
The product $\be_{mi}^{(s)}\neq0$ and is well defined.
\endproclaim

\demo{Proof}
Any factor in the product $\be_{mi}^{(s)}$ corresponding to $r\neq s$ 
is non-zero and well defined since $h^{(r)}-h^{(s)}\notin\ZZ$ then. Since
$\La^{(s)},\La_\tz^{(s)}\in\S_{\la^{(s)},\mu^{(s)}}$
we have the inequalities
$$
\alignat2
&
\ka_{M^{(s)}+m+1,j}^{(s)}
\geqslant
\la_{M^{(s)}+m+1,j}^{(s)}
\geqslant
\la_{M^{(s)}+m,j}^{(s)}
\geqslant
\la_{M^{(s)}+m,i}^{(s)}
&&
\qquad\text{if}\quad j\leqslant i\ts,
\\
&
\ka_{M^{(s)}+m-1,j}^{(s)}
\geqslant
\la_{M^{(s)}+m-1,j}^{(s)}
\geqslant
\la_{M^{(s)}+m,j+1}^{(s)}
\geqslant
\la_{M^{(s)}+m,i}^{(s)}
&&
\qquad\text{if}\quad j< i\ts.
\endalignat
$$
These inequalities show that
any factor in the product $\be_{mi}^{(s)}$ corresponding to $r=s$ 
is also non-zero and well defined
\enddemos

\proclaim{Theorem 3.8}
We have 
$$
b_m\bigl(\nu_{mi}^{(s)}\bigr)\cdot\xi=
\cases
\ts\ts\be_{mi}^{(s)}\cdot\xi_-\quad
&\text{if}\quad\La_\tm^{(s)}\in\S_{\la^{(s)},\mu^{(s)}}\ts;
\\
\ \ts\quad0\quad
&\text{otherwise}.
\endcases
$$
\endproclaim

\demo{Proof}
We will use again some of the arguments which appeared in the proofs of
Theorem
3.2 and Proposition 3.3. Consider the vector
$b_m\bigl(\nu_{mi}^{(s)}\bigr)\cdot\xi\in W$.
Let us check first that for any $l=1\lc N$ we have 
$$
a_l(u)\ts b_m(\nu_{mi}^{(s)})\cdot\xi=
\varpi_{l,\La^{(1)}\lc\La_\tm^{(s)}\lc\La^{(n)}}(u)\ts
b_m(\nu_{mi}^{(s)})\cdot\xi.
\Tag{3.8}
$$
Indeed, if $l\neq m$ then
$$
\varpi_{l,\La^{(1)}\lc\La_\tm^{(s)}\lc\La^{(n)}}(u)=
\varpi_{l,\La^{(1)}\lc\La^{(s)}\lc\La^{(n)}}(u).
$$
On the other hand, by the relation \(1.7) and by Theorem 3.2 we then have
$$
a_l(u)\ts b_m(\nu_{mi}^{(s)})\cdot\xi=
b_m(\nu_{mi}^{(s)})\ts a_l(u)\cdot\xi=
\varpi_{l,\La^{(1)}\lc\La^{(s)}\lc\La^{(n)}}(u)\ts b_m(\nu_{mi}^{(s)})
\cdot\xi\ts.
$$
Suppose that $l=m\ts$; then by the definition of $\La_\tm^{(s)}$ we have
$$
\varpi_{m,\La^{(1)}\lc\La_\tm^{(s)}\lc\La^{(n)}}(u)=
\frac{u-\nu_{mj}^{(s)}-1}{u-\nu_{mj}^{(s)}}\ts\ts
\varpi_{m,\La^{(1)}\lc\La^{(s)}\lc\La^{(n)}}(u).
$$
Since due to Theorem 3.2
$$
a_m(\nu_{mi}^{(s)})\cdot\xi=
\varpi_{m,\La^{(1)}\lc\La^{(s)}\lc\La^{(n)}}(\nu_{mi}^{(s)})\cdot\xi=0,
\Tag{3.85}
$$
by the relation
\(1.10) and again by Theorem 3.2 we get the equalities
$$
\align
a_m(u)\ts b_m(\nu_{mi}^{(s)})\cdot\xi
&=
\frac{u-\nu_{mi}^{(s)}-1}{u-\nu_{mi}^{(s)}}\ts\ts
b_m(\nu_{mi}^{(s)})\ts a_m(u)\cdot\xi
\\
&=
\varpi_{m,\La^{(1)}\lc\La_\tm^{(s)}\lc\La^{(n)}}(u)
\ts b_m(\nu_{mi}^{(s)})\cdot\xi\ts.
\endalign
$$

Since
$\La^{(s)}\in\S_{\la^{(s)},\mu^{(s)}}$
the array $\La_\tm^{(s)}$ can be uniquely restored from the collection of
the polynomials
$$
\varpi_{l,\La^{(1)}\lc\La_\tm^{(s)}\lc\La^{(n)}}(u);
\qquad
l=1\lc N.
$$
Thus if $\La_\tm^{(s)}\notin\S_{\la^{(s)},\mu^{(s)}}$ then 
$b_m(\nu_{mi}^{(s)})\cdot\xi=0$ by Theorem 2.9
as we have claimed.

Assume that $\La_\tm^{(s)}\in\S_{\la^{(s)},\mu^{(s)}}$.
Then due to Theorem 2.9 and Proposition 3.6 the equalities \(3.8) imply that
$$
b_m(\nu_{mi}^{(s)})\cdot\xi=\be\cdot\xi_- 
\Tag{3.9}
$$
for some $\be\in\CC$.
We will prove that $\be=\be_{mi}^{(s)}$ here.

Let us compare the action of the element
$c_m(\nu_{mi}^{(s)}+1)$ on the both sides of the equality \(3.9).
On the right hand side wy Theorem 3.5 we have
$$
c_m(\nu_{mi}^{(s)}+1)\ts\be\cdot\xi=-\be\ts\ga\cdot\xi
$$
where $\ga$ is analogue of $\ga_{mi}^{(s)}$ for the collection of
schemes
$\La^{(1)}\lc\La_\tm^{(s)}\lc\La^{(n)}$ instead of $\La^{(1)}\lc
\La^{(s)}\lc\La^{(n)}$.
In particular, here $\ga\neq0$ 
by Lemma 3.4. On the left hand side by applying
Proposition 1.10 and using the equality \(3.85)
we get
$$
\align
c_m(\nu_{mi}^{(s)}+1)\ts b_m(\nu_{mi}^{(s)})\cdot\xi
=
-a_{m+1}(\nu_{mi}^{(s)}+1)\ts a_{m-1}(\nu_{mi}^{(s)})
&\cdot\xi
\\
=-\ts
\varpi_{m+1,\La^{(1)}\lc\La^{(s)}\lc\La^{(n)}}(\nu_{mi}^{(s)}+1)\ts
\varpi_{m-1,\La^{(1)}\lc\La^{(s)}\lc\La^{(n)}}(\nu_{mi}^{(s)})
&\cdot\xi\ts.
\endalign
$$
Since $\xi\neq0$ by Proposition 3.6, we finally obtain that
$$
\be=
\varpi_{m+1,\La^{(1)}\lc\La^{(s)}\lc\La^{(n)}}(\nu_{mi}^{(s)}+1)\ts
\varpi_{m-1,\La^{(1)}\lc\La^{(s)}\lc\La^{(n)}}(\nu_{mi}^{(s)})\ts\ga^{-1}
=\be_{mi}^{(s)}\ts.
$$
Thus we have proved Theorem 3.8
\enddemos

\nt
Due to Proposition 3.6 we can choose $\xi_{\La^{(1)}\lc\La^{(n)}}=\xi$.
Theorems 3.2 and 3.5\ts,\ts3.8 then
completely describe the action of the Drinfeld generators of the algebra
$\YN$
in the module $W$. 
Since the vector $\xi_{\La^{(1)}\lc\La^{(n)}}$ is uniquely determined by
the equalities
\(2.3), we obtain the following corollary to Lemmas 3.4\ts,\ts 3.7.
\nopagebreak
\proclaim{Corollary 3.9}
$\YN$-module $W$ is irreducible
if $h^{(r)}-h^{(s)}\notin\ZZ$~for~all
\nolinebreak
$r~\neq~s$.
\endproclaim

\bigskip\bigskip
\Section{Classification theorem}

\nt
Let $V$ be any irreducible finite-dimensional
module over the algebra $\YN$. Let $P_1(u)\lc P_{N-1}(u)$ be the
Drinfeld polynomials corresponding to $V$. For each $m=1\lc N-1$
consider the collection of zeroes of the polynomial $P_m(-u)$
$$
\bigl(\ts z_{mi}\ts|\ts i=1\lc\deg P_m\ts\bigr).
$$

\proclaim{Theorem 4.1}
The following three conditions are equivalent.
\smallskip
\itemitem{a)}
For all $m\geqslant l$ we have $z_{li}-z_{mj}\neq 0\lc m-l$ unless
$(m,i)=(l,j)$.
\itemitem{b)}
Up to some automorphism $\of$ of $\YN$
the module $V$ has the form \(2.29) where
$h^{(r)}-h^{(s)}\notin\ZZ$ for all $r\neq s$.
\itemitem{c)}
The action in the module $V$ of the subalgebra $\AN$ is semisimple. 
\endproclaim

\demo{Proof}
Due to Theorem 2.9 the condition (b) implies (c). We will demonstrate that
(a) implies (b) and that (c) implies (a). 
Suppose that the condition (a) is
satisfied. In particular, then we have 
$u_{mi}\neq u_{lj}$ unless $(m,i)=(l,j)$.
Consider the partition of the set of all zeroes
$$
\Z=
\bigl\{\ts z_{mi}\ts|\ts m=1\lc N-1;\ i=1\lc\deg P_m\ts\bigr\}
=
\Z^{(1)}\sqcup\ldots\sqcup\Z^{(n)}
$$
where $\Z^{(1)}\lc\Z^{(n)}$ are the maximal subsets in $\Z$ such that
$$
z\ts,w\in\Z^{(s)}\ts\ \Rightarrow\ts\ z-w\in\ZZ\ts;\qquad s=1\lc n.
$$
Furthermore for each $s=1\lc n$ introduce the partition 
$$
\Z^{(s)}=\Z^{(s)}_{\ts 1}\sqcup\ldots\sqcup\Z^{(n)}_{\ts N-1}
$$
into the subsets of zeroes of the polynomials $P_1(-u)\lc P_{N-1}(-u)$
respectively. 
We will suppose that $\Z\neq\varnothing$; otherwise the module $V$ is
one\ts-dimensional and there we have $T_{ij}(u)\mapsto f(u)\cdot\de_{ij}$
for
certain series $f(u)\in 1+u^{-1}\ts\CC[[u^{-1}]]$. For each $s=1\lc n$
we will
point out a parameter $h^{(s)}\in\CC$ and a skew Young diagram
$\la^{(s)}/\mu^{(s)}$
whose columns are of the height at most $N-1$ such that
$$
P_{m,\la^{(s)}/\mu^{(s)}}\bigl(u+h^{(s)}\bigr)=
\prod_{z\ts\in\ts\Z^{(s)}_{\ts m}}\ts(u+z);
\qquad
m=1\lc N-1.
\Tag{4.-1}
$$
Then by Corollary 2.13 the module $V$ will have the same Drinfeld
polynomials as the
module \(2.29). Therefore the condition (b) will be also satisfied.

Let any $s\in\{\ts 1\lc n\ts\}$ be fixed such that
$\Z^{(s)}\neq\varnothing$. 
Denote $\#\ts\Z^{(s)}=p$. Let
$z_1\lc z_p$ be the elements of the set $\Z^{(s)}\ts$; we will assume that
$z_k-z_{k-1}<0$ for all possible indices $k$. Put $h^{(s)}=z_1$.
Suppose that $z_1\in\Z^{(s)}_{\ts q}\ts$. Then assign to each $k=1\lc p$
the column
$$
\bigl\{\ts(i,j)\in\ZZ^2\ts|\ts\ts
j=q-k+1,
\ 
z_k-z_1\leqslant j-i<z_k-z_1+m
\ts\bigr\}
\Tag{4.0}
$$
where $z_k\in\Z^{(s)}_{\ts m}$. Note that here for $k=1$ we have
$i\geqslant1$.
By our assumption for any $k>1$ we have the inequality
$$
z_k-z_1<z_{k-1}-z_1\ts.
\Tag{4.1}
$$
Moreover, if $z_{k-1}\in\Z^{(s)}_{\ts l}$ we also have the inequality
$$
z_k-z_1+m<z_{k-1}-z_1+l\ts.
\Tag{4.2}
$$
If $m\leqslant l$ it follows from \(4.1). If $m>l$ then
$z_{k-1}-z_k\neq 1\lc m-l$
by the condition (a) so that \(4.2) again follows from \(4.1).
Therefore the collection
of all columns \(4.0) is a skew Young diagram. Denote this diagram by
$\la^{(s)}/\mu^{(s)}\ts$; the equalities \(4.-1) then evidently hold.
Thus we have proved that (a) implies (b).

It remains to prove that the condition (c) implies (a). We will consider
first the
case $N=2$. Denote by
$V_k(h)$ the $\operatorname{Y}(\gl_2)$-module  
$V_\lm(h)$ with $\la=(k,0)$ and $\mu=\varnothing$.
Then up to some authomorphism $\of$ every irreducible finite-dimensional  
$\operatorname{Y}(\gl_2)$-module has the form
$$
W=V_{k^{(1)}}\bigl(h^{(1)}\bigr)\ot\ldots\ot V_{k^{(n)}}\bigl(h^{(n)}\bigr)
$$
for certain $k^{(1)}\lc k^{(n)}\in\{\ts1,2\ldots\ts\}$
and $h^{(1)}\lc h^{(n)}\in\CC$.
Moreover, if we denote
$$
\X^{(s)}=\bigl\{\ts h^{(s)},h^{(s)}+1\lc h^{(s)}+k^{(s)}\ts\bigr\}\ts;
\qquad
s=1\lc n
$$
then the condition $\X^{(r)}\cap\X^{(s)}\neq\varnothing$ implies that
either $\X^{(r)}\subset\X^{(s)}$ or $\X^{(r)}\supset\nomathbreak\X^{(s)}$. 
These two facts are contained
in [T2]; see also [CP1], Theorems 4.1 and 4.11. 
Furthermore, by Proposition 2.8 for any choice of the elements
$x^{(s)}\in\X^{(s)}$ with $s=1\lc n$
there exists a non-zero vector $\xi\in W$ such that
$$
A_1(u)\cdot\xi=\xi\cdot\ts
\prod_{s=1}^n\ts\ts\frac{u+x^{(s)}}{u+h^{(s)}}\ts,
\quad
A_2(u)\cdot\xi=\xi\cdot\ts\prod_{s=1}^n\ts\ts
\frac{u+k^{(s)}+h^{(s)}}{u+h^{(s)}}\ts.
\Tag{4.25}
$$

Put $\rho_1(u)=(u+h^{(1)})\ldots(u+h^{(n)})$.
Consider the formal Laurent series in $u^{-1}$
$$
\rho_1(u)\ts A_1(u),
\quad
\rho_1(u)\ts B_1(u),
\quad
\rho_1(u)\ts C_1(u),
\quad
\rho_1(u)\ts D_1(u).
\Tag{4.3}
$$
By the definition of the comultiplication \(1.135) 
the action in $W$
of every coefficient of these series at $u^{-1},u^{-2},\ldots$
vanishes. Denote by $a(u),b(u),c(u),d(u)$
the polynomials in $u$ with the coefficients in $\End(W)$ corresponding
to the
series \(4.3).

Suppose that the Drinfeld polynomial $P(u)$ of the
$\operatorname{Y}(\gl_2)$-module $W$ has multiple zeroes.
We have to show that the action in $W$ of the subalgebra 
$\operatorname{A}(\gl_2)$ is then not semisimple.
Due to Corollary 2.13 we have
$$
P(u)\ts=\ts\prod_{s=1}^n\ts\prod_{j=1}^{k^{(s)}}\ts
\bigl(u+j-1+h^{(s)}\bigr).
$$
Let
$z$ be the zero of the polynomial $P(-u)$ of the maximal multiplicity.
For some $r\in\{\ts1\lc n\ts\}$
we have $z\in\X^{(r)}$. We will assume that the set $\X^{(r)}$ is minimal
amongst
those of $\X^{(1)}\lc\X^{(n)}$ which contain $z$. Then
the condition $\X^{(r)}\cap\X^{(s)}\neq\varnothing$ implies that
$\X^{(r)}\subset\X^{(s)}$. In particular, for $x=k^{(r)}+h^{(r)}$ we have
$$
x\neq h^{(s)},\ts k^{(s)}+h^{(s)}+1\ts;
\qquad
s=1\lc n.
\Tag{4.325}
$$
Now set
$$
x^{(s)}=
\cases
x&\quad\text{if}\ x\in\X^{(s)};
\\
h^{(s)}&\quad\text{otherwise}
\endcases
$$
and consider a non-zero vector $\xi\in W$ satisfying the corresponding
equalities \(4.25).
Put $\al(u)=(u+x^{(1)})\ldots(u+x^{(n)})$. Then we have 
$a(u)\cdot\xi=\al(u)\cdot\xi$. In particular,
$$
a(-x)\cdot\xi=\al(-x)\cdot\xi=0\ts,
\qquad
\frac{\deri a}{\deri u}(-x)\cdot\xi=
\frac{\deri \al}{\deri u}(-x)\cdot\xi=0\ts.
\Tag{4.35}
$$
Introduce the vectors in $W$
$$
\xi^{\ts\prime}=b(-x)\cdot\xi\ts,
\qquad
\xi^{\ts\prime\prime}=\frac{\deri b}{\deri u}(-x)\cdot\xi\ts.
$$
We shall prove that 
$$
\align
a(u)\cdot\xi^{\ts\prime}
&=
\al(u)\ts\frac{u+x-1}{u+x}\cdot\xi^{\ts\prime}\ts;
\Tag{4.4}
\\
a(u)\cdot\xi^{\ts\prime\prime}
&=
\al(u)\ts\frac{u+x-1}{u+x}\cdot\xi^{\ts\prime\prime}
-
\frac{\al(u)}{(u+x)^2}\cdot\xi^{\ts\prime}\ts.
\Tag{4.5}
\endalign
$$
and that $\xi^{\ts\prime}\neq0$. Then the equalities \(4.4),\ts\(4.5)
will imply
that the action in $W$ of the subalgebra in
$\operatorname{Y}(\gl_2)$ generated by the coefficients of $A_1(u)$
is not semisimple.

Due to \(4.35) by the relation \(1.10) we have
$$
\align
a(u)\cdot\xi^{\ts\prime}
&
=a(u)\ts b(-x)\cdot\xi
=\frac{u+x-1}{u+x}\ts b(-x)\ts a(u)\cdot\xi
\\
&
=\al(u)\ts\frac{u+x-1}{u+x}\ts b(-x)\cdot\xi
=\al(u)\ts\frac{u+x-1}{u+x}\cdot\xi^{\ts\prime}\ts\ts;
\\
\vspace{6pt}
a(u)\cdot\xi^{\ts\prime\prime}
&
=a(u)\ts
\frac{\deri b}{\deri u}(-x)\cdot\xi
\\
&
=\frac{u+x-1}{u+x}\ts
\frac{\deri b}{\deri u}(-x)\ts a(u)\cdot\xi
+
\frac1{u+x}\ts[\ts a(u),b(-x)]\cdot\xi
\\
&
=\al(u)\ts\frac{u+x-1}{u+x}\ts
\frac{\deri b}{\deri u}(-x)\cdot\xi
-
\frac{\al(u)}{(u+x)^2}\ts b(-x)\cdot\xi
\\
&
=
\al(u)\ts\frac{u+x-1}{u+x}\cdot\xi^{\ts\prime\prime}
-
\frac{\al(u)}{(u+x)^2}\cdot\xi^{\ts\prime}\ts.
\endalign
$$
Thus we have proved \(4.4) and \(4.5).
Furthermore, due to the second equality in \(4.25) by the relation
\(1.12) we have
$$
\align
c(1-x)\cdot\xi^{\ts\prime}
&
=c(1-x)\ts b(-x)\cdot\xi
=\rho_1(1-x)\ts\rho_1(-x)\ts A_2(1-x)\cdot\xi
\\
&
=\xi\cdot\ts\prod_{s=1}^n\ts\ts
\bigl(k^{(s)}+h^{(s)}+1-x\bigr)\ts\bigl(h^{(s)}-x\bigr)\ts.
\endalign
\nopagebreak
$$
Every factor in the above product differs from zero due to \(4.325) so that
$\xi^{\ts\prime}\neq0$.

Thus if the action of  
$\operatorname{A}(\gl_2)$ in an irreducible finite-dimensional
$\operatorname{Y}(\gl_2)$-module is semisimple then the corresponding
Drinfeld
polynomial has no multiple zeroes. Hence (c) implies (a) for $N=2$.
Let us now prove that (c) implies (a) for $N\geqslant3$.

Let $\xi\in V$ be a singular vector.
Suppose that the action in the module $V$ of the subalgebra $\AN$ is
semisimple.
Note that due to \(1.655555) 
the action of $\AN$ in the $\YN$-module $\overline{V}$ 
is then also semisimple.

Let the indices $m,l\in\{1\lc N-1\}$ such that $m\geqslant l$ be fixed. 
Consider the embedding $\ph^{\vp1}_{l,m+1}$ of the algebra
$\operatorname{Y}(\gl_2)$
into $\YN$. By Proposition 1.14~the~action of
$\operatorname{A}(\gl_2)$
in $V$ corresponding to this embedding is semisimple.
The action of $\operatorname{A}(\gl_2)$
in the irreducible $\operatorname{Y}(\gl_2)$-subquotient of $V$ generated
by the vector $\xi$ is then also semisimple.
But due to Proposition 2.10 in the module $V$ we have
$$
\gather
\frac{T_{m+1,m+1}(u-1)}{T_{l,l}(u-1)}\cdot\xi=
\\
\frac{T_{m+1,m+1}(u-1)}{T_{m,m}(u-1)}
\cdot
\frac{T_{m,m}(u-1)}{T_{m-1,m-1}(u-1)}
\cdot
\ \ldots\ 
\cdot
\frac{T_{l+1,l+1}(u-1)}{T_{l,l}(u-1)}
\cdot\xi=
\\
\frac{P_m(u+m-2)}{P_m(u+m-1)}
\cdot
\frac{P_{m-1}(u+m-3)}{P_{m-1}(u+m-2)}
\cdot
\ \ldots\ 
\cdot
\frac{P_l(u+l-2)}{P_l(u+l-1)}
\cdot\xi
\endgather
$$
so that the Drinfeld polynomial of the above introduced subquotient equals
$$
P_m(u+m-1)\ts P_{m-1}(u+m-2)\ts\ldots\ts P_l(u+l-1).
$$
As we have already proved all zeroes of this polynomial are simple. 
Therefore
$$
z_{li}-z_{mj}\neq m-l\quad\text{unless}\quad(m,i)=(l,j).
\Tag{4.6}
$$
By Proposition 2.11 the Drinfeld polynomial of the irreducible
$\operatorname{Y}(\gl_2)$-subquotient of the module $\overline V$ generated
by the vector $\xi$ then equals
$$
\align
\overline{P}_m(u+m-1)\ts\overline{P}_{m-1}(u+m-2)\ts
&\ldots\ts\overline{P}_l(u+l-1)
=
\\
P_{N-m}(u+N-1)\ts P_{N-m+1}(u+N-1)\ts&\ldots\ts P_{N-l}(u+N-1).
\endalign
$$
All zeroes of this polynomial are also simple. Replacing here $N-m$ and
$N-l$~by
$l$~and~$m$ respectively we obtain that
the polynomials $P_l(u)$ and $P_m(u)$~of~the module $V$ have no common 
zeroes.

Let us now assume that $m>l$. Fix an arbitrary $k\in\{1\lc m-l\}$ and put
$$
\bk=(l,l+1\lc l+k,m+1).
$$
Consider the embedding $\ph^{\vp1}_\bk$ of the algebra
$\operatorname{Y}(\gl_{k+2})$ into $\YN$.
By Proposition 1.14~the~action of
$\operatorname{A}(\gl_{k+2})$
in $V$ corresponding to this embedding is semisimple.
The action of $\operatorname{A}(\gl_{k+2})$
in the irreducible $\operatorname{Y}(\gl_{k+2})$-subquotient of $V$
generated
by the vector $\xi$ is then also semisimple.
Due to Proposition 2.10 the first and the last Drinfeld polynomials of
this
subquotient are respectively
$P_l(u+l-1)$ and
$$
P_m(u+m-k-1)\ts P_{m-1}(u+m-k-2)\ts\ldots\ts P_{l+k}(u+l-1).
$$
As we have already proved, these two polynomials have no common zeroes.
Thus
$$
z_{li}-z_{mj}\neq m-l-k\quad\text{unless}\quad(m,i)=(l,j)\ts;
\qquad
k=1\lc m-l.
\nopagebreak
\Tag{4.7}
$$
The statements \(4.6) and \(4.7) constitute the condition a). We have
completed the
proof of Theorem 4.1 
\enddemos

\bigskip\bigskip
\section{References}

\itemitem{[C]}
{\it I. Cherednik},
{A new interpretation of Gelfand-Zetlin bases},
{Duke Math. J.}
{\bf 54}
(1987),
563--577.

\itemitem{[CP1]}
{\it V. Chari and A. Pressley},
{Yangians and $R$-matrices},
{L'Enseign. Math.}
{\bf 36}
(1990),
267--302.

\itemitem{[CP2]}
{\it V. Chari and A. Pressley},
{Fundamental representations of Yangians and
singularities of $R$-matrices},
{J. Reine Angew. Math.}
{\bf 417}
(1991),
87--128.

\itemitem{[D1]}
{\it V. Drinfeld},
{Hopf algebras and the quantum Yang--Baxter equation},
{\text{Soviet} Math. Dokl.}
{\bf 32}
(1985),
254--258.

\itemitem{[D2]}
{\it V. Drinfeld},
{A new realization of Yangians and quantized affine algebras},
{Soviet Math. Dokl.}
{\bf 36}
(1988),
212--216.

\itemitem{[D3]}
{\it V. G. Drinfeld},
Quantum Groups,
in {\lq\lq\ts 
Proceedings of the Internation\-al Con\-gress of Mathematicians
\ts\rq\rq},
Amer. Math. Soc.,
Providence,
1987,
798--820.

\itemitem{[GRV]}
{\it V. Ginzburg, N. Reshetikhin and E. Vasserot},
{Quantum groups and flag varieties}, 
Contemp. Math.
{\bf 175}
(1994),
101-130.

\itemitem{[GZ]}
{\it I. Gelfand and M. Zetlin},
{Finite-dimensional representations of the unimodular group},
{Dokl. Akad. Nauk SSSR}
{\bf 71}
(1950),
825--828.

\itemitem{[M]}
{\it A. Molev},
{Gelfand-Tsetlin bases for representations of Yangians},
{Lett. Math. Phys.}
{\bf 30}
(1994),
53--60.

\itemitem{[MNO]}
{\it A. Molev, M. Nazarov and G. Olshanski},
{Yangians and classical Lie algebras}, 
{Russian Math. Surveys}
{\bf 51}
(1996),
205--282.

\itemitem{[NT]}
{\it M. Nazarov and V. Tarasov},
{Yangians and Gelfand-Zetlin bases},
{Publ. RIMS}
{\bf 30}
(1994),
459--478.

\itemitem{[O]}
{\it G. Olshanski},
{Representations of infinite-dimensional classical groups,
limits of enveloping algebras, and Yangians},
{Adv. Soviet Math.}
{\bf 2},
1991,
1--66.

\itemitem{[T1]}
{\it V. Tarasov},
{Structure of quantum {\it L\/}-operators for the
{\it R\/}-matrix of the {\it XXZ\/}-model},
{Theor. Math. Phys.}
{\bf 61}
(1984),
1065--1071.

\itemitem{[T2]}
{\it V. Tarasov},
{Irreducible monodromy matrices for the {\it R\/}-matrix of the
{\it XXZ\/}-model and lattice local quantum Hamiltonians},
{Theor. Math. Phys.}
{\bf 63}
(1985),
440--454.
\bigskip
\centerline{\hbox to 4cm{\hrulefill}}
\bigskip
\centerline{Mathematics Department, University of York,
Heslington, York YO1 5DD, England}
\centerline{e-mail: mln1\@ york.ac.uk}
\medskip
\centerline{Steklov Mathematical Institute, Fontanka 27,
St.\,Petersburg 191011, Russia}
\centerline{e-mail: vt\@ pdmi.ras.ru}

\bye